%% file: main.tex
\documentclass[11pt]{article}

\usepackage{amsthm}
\usepackage{amsmath, bm}
\usepackage{amssymb}
\usepackage{natbib}
\usepackage{geometry}
\usepackage{subfig}
\usepackage{graphicx}
\usepackage{tikz}
\usepackage{multirow}
\usepackage{listings}
\usepackage{xcolor}
\usepackage{authblk}
\usepackage{float}
\usepackage{url}
\usepackage{calc}
\usepackage{stmaryrd}

\usepackage{subfiles}
\usepackage{indentfirst}
\usepackage{fullpage}

\newtheorem{Def}{Definition}[section]
\newtheorem{Thm}{Theorem}[section]

\newtheorem{theorem}{Theorem}[section]
\newtheorem{lemma}[theorem]{Lemma}

\newcommand{\diag}[1]{{\rm diag}\LRp{#1}}
\newcommand{\td}[2]{\frac{{\rm d}#1}{{\rm d}{\rm #2}}}
\newcommand{\pd}[2]{\frac{\partial#1}{\partial#2}}
\newcommand{\nor}[1]{\left\| #1 \right\|}
\newcommand{\LRp}[1]{\left( #1 \right)}
\newcommand{\LRs}[1]{\left[ #1 \right]}

\newcommand{\LRc}[1]{\left\{ #1 \right\}}

\newcommand{\avg}[1] {\ensuremath{\LRc{\!\{#1\}\!}}}

\newcommand{\fnt}[1]{\bm{\mathsf{ #1}}}

\title{Entropy stable discontinuous Galerkin methods for nonlinear conservation laws on networks and multi-dimensional domains}
\author[1]{Xinhui Wu}
\author[1]{Jesse Chan}
\affil[1]{Department of Computational and Applied Mathematics, Rice University}
\date{}

\begin{document}
\maketitle

\begin{abstract}
We present a high-order entropy stable discontinuous Galerkin (ESDG) method for nonlinear conservation laws on both multi-dimensional domains and on networks constructed from one-dimensional domains. These methods utilize treatments of multi-dimensional interfaces and network junctions which retain entropy stability when coupling together entropy stable discretizations. Numerical experiments verify the stability of the proposed schemes, and comparisons with fully 2D implementations demonstrate the accuracy of each type of junction treatment.
\end{abstract}

\subfile{section1}
\subfile{section2}
\subfile{section3}
\subfile{section4}
\subfile{section5}
\subfile{section6}

\section*{Acknowledgement}
Authors Xinhui Wu and Jesse Chan gratefully acknowledge support from the National Science Foundation under awards DMS-1719818, DMS-1712639, and DMS-CAREER-1943186.



\subfile{appendixA}

\bibliographystyle{unsrt}
\bibliography{main.bib}

\end{document}

%% file: section1.tex
\section{Introduction}
\label{sec:Introduction}
\begin{subequations}
\numberwithin{equation}{section}

There is an increasing interest in the mathematical modeling of physical systems posed on spatial domains with network-like structures. In this situation, the 2D (or 3D) partial differential equations (PDEs) that describe the physics of the system can be well-approximated by 1D PDEs. Simulations over network-like domains can then be performed by coupling together 1D subdomains at junctions. This reduction from fully 2D or 3D simulations to simulations over 1D domains both simplifies the construction of meshes for network-like domains and reduces computational cost \cite{bellamoli2018numerical, sherwin2003computational, neupane2015discontinuous, west2017multidimensional, borsche2014ader}. 

Applications of network models include the simulation of gas flow in pipelines \cite{banda2006coupling, brouwer2011gas, reigstad2015coupling}, water flow through channels \cite{kesserwani2008simulation, akan1981diffusion, zhang2005simulation}, and blood flow in the human cardiovascular system \cite{sherwin2003computational, taylor1998finite, quarteroni2000computational, muller2015high}. In gas networks, flow is governed by equations derived from the compressible Euler equations. In water flow through open channels, the physics are described by the shallow water equations. 
In blood vessels, the system is governed by a system of equations which closely resembles the shallow water equations \cite{sherwin2003computational}. 
For each of these examples, however, the systems involved are nonlinear conservation laws. To demonstrate the idea, we present the fully 2D discretization of a river with turning channel and its 1D-2D model in Figure \ref{fig:demo}.

\begin{figure}
\centering
\subfloat[]{\includegraphics[width=.4\textwidth]{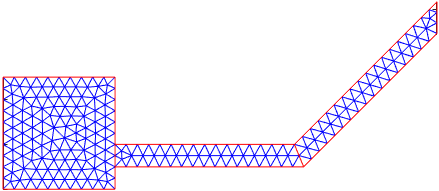}}
\subfloat[]{\includegraphics[width=.4\textwidth]{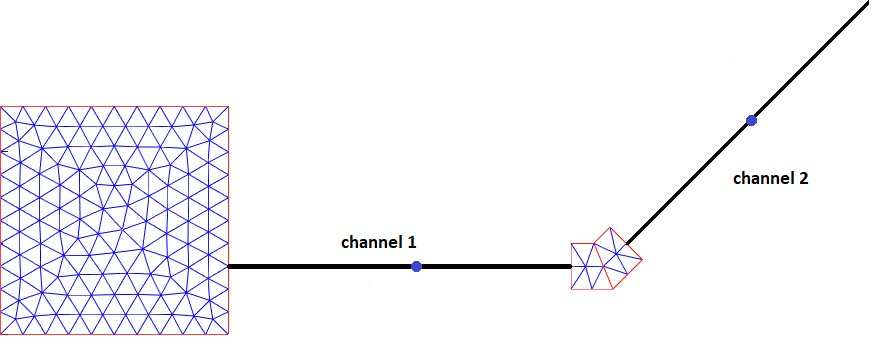}}
\caption{Lake and turning channel meshes. Fully 2D mesh (a), 1D-2D (b)}
\label{fig:demo}
\end{figure}

In this paper, we present new numerical treatments of one-dimensional junctions and interfaces between multi-dimensional (1D and 2D) domains. When combined with entropy stable discretizations, these interface treatments preserve an entropy inequality at the semi-discrete level. Individual 1D and 2D domains are discretized using entropy stable discontinuous Galerkin (DG) methods, which provide geometric flexibility for 2D subdomains through the use of unstructured meshes \cite{neupane2015discontinuous, chen2017entropy}. It is also straightforward to increase the order of accuracy of DG methods by increasing the degree of the polynomial approximation \cite{chan2018discretely}. Finally, DG methods have been shown to be highly parallelizable and amenable to high performance implementations on graphics processing units (GPUs) \cite{chan2016gpu, chan2017bbdg, wintermeyer2018entropy}. 

There exist several different entropy stable DG schemes for the shallow water and compressible Euler equations \cite{chen2017entropy, wintermeyer2018entropy, chan2018discretely}. Traditional entropy stable DG formulations have relied on finite difference summation-by-parts (SBP) operators, which are constructed using carefully designed quadrature rules which contain boundary points while satisfying certain accuracy conditions \cite{chen2017entropy}. In this work, we utilize entropy stable ``modal'' DG formulations as introduced by \cite{chan2018discretely, chan2019skew}. However, our focus is on the construction of entropy stable treatments of multi-dimensional interfaces and junctions, and are applicable to both ``modal'' and traditional entropy stable SBP-DG schemes. We also note that, while we consider only affine triangular meshes, the approaches presented here can be extend to curved meshes of either triangles or quadrilaterals so long as the 1D-2D interface is a straight line.



This paper is organized as follows: we begin by reviewing nonlinear conservation laws and entropy theory in Section \ref{sec:2}. In Section \ref{sec:3}, we review the construction of entropy stable DG methods. In Section \ref{sec:4}, we construct entropy stable couplings for multidimensional interfaces and junctions. We present numerical experiments which validate theoretical properties of our schemes in Section \ref{sec:numerical_results}. Finally, we conclude in Section \ref{sec:Conclusion} with a summary of the main results.

\end{subequations}

%% file: section2.tex
\section{Entropy inequalities for nonlinear conservation laws}
\label{sec:2}
\begin{subequations}
\numberwithin{equation}{section}

In this paper, we focus on numerical methods for system of nonlinear conservation laws. The general form of such a system in $d$ dimensions is
\begin{align}
\frac{\partial \bm{u}}{\partial t} + \sum_{i=1}^d\frac{\partial \bm{f}_i(\bm{u})}{\partial x_i} = 0, \hspace{1cm} (\bm{x},t) \in \mathbb{R}^d\times[0,\infty],
\label{eq:nonlinear_general}
\end{align}
where $\bm{x}\in\Omega\subseteq\mathbb{R}^d$, $\bm{u}(\bm{x}, t)  =[u_1, ...u_n]^T$ denote the conservative variables and $\bm{f}_i(\bm{u})$ denote the flux functions. We assume that there exists a convex scalar entropy function $S(\bm{u})$ such that
\begin{align}
S''(\bm{u})\bm{A}(\bm{u}) = \LRp{S''(\bm{u})\bm{A}(\bm{u})}^T, \hspace{1cm} \bm{A}(\bm{u})_{ij} = \LRp{\frac{\partial \bm{f}(\bm{u})}{\partial u_j}}_i,
\label{eq:entropy_function}
\end{align}
where $\bm{A}(\bm{u})$ is the flux Jacobian matrix. 

We next define the entropy variables $\bm{v} = S'(\bm{u})$. For values of $\bm{u}$ over which the entropy function is convex, the mapping between conservative variables and entropy variables is invertible. Then, it can be shown \cite{mock1980systems} that there exists entropy flux functions $F_i(\bm{u})$ and entropy potentials $\psi_i(\bm{u})$ such that
\begin{align}
\bm{v}(\bm{u})^T\frac{\partial \bm{f}_i(\bm(u))}{\partial \bm{u}} = \frac{\partial F_i(\bm{u})^T}{\partial \bm{u}}, \hspace{.5cm}
\psi_i(\bm{u}) = \bm{v}(\bm{u})^T\bm{f}_i(\bm{u}) - F_i(\bm{u}), \hspace{0.5cm} \psi'_i(\bm{u}) = \bm{f}_i(\bm{u}).
\end{align}
In regions where the solution is smooth, an entropy equality can be derived for $\bm{u}$ by multiplying the conservation law by $\bm{v}^T$ and integrating over the domain. Then, using the chain rule and definition of the entropy flux, we have the following statement of entropy conservation on domain $\Omega$
\begin{align}
\int_{\Omega}\frac{\partial S(\bm{u})}{\partial t} = \bm{0}.
\label{eq:ec_bottom0}
\end{align} 

In this paper, we will construct a high order DG scheme for multi-dimensional and network domains that is entropy conservative at the semi-discrete level. By adding appropriate entropy dissipation terms at inter-element interfaces, multi-dimensional interfaces, and junctions, these entropy conservative schemes can be made entropy stable. In this work, we focus specifically on 1D junction treatments and coupling between 1D-2D domains for the shallow water equations (SWE) and the compressible Euler equations. 

\subsection{Shallow water equations in one and two dimensions} 

We begin with introducing the two-dimensional shallow water equations 
\begin{align}
\frac{\partial }{\partial t}\begin{bmatrix}
h\\
hu\\
hv
\end{bmatrix}+
\frac{\partial }{\partial x}\begin{bmatrix}
hu\\
hu^2+gh^2/2\\
huv
\end{bmatrix}+
\frac{\partial }{\partial y}\begin{bmatrix}
hv\\
huv\\
hv^2+gh^2/2
\end{bmatrix} = 0.
\label{eq:SWE_2D}
\end{align}
Here, $h$ denotes the water height as measured from the lake or channel bottom. The velocity in the $x$ direction is denoted by $u$ and the velocity in the $y$ direction is denoted by $v$. The gravitational constant is denoted by $g$.  In this example, we have the conservative variable $\bm{u} = [h,hu,hv]^T$ and the flux functions are $\bm{f}_1 = [hu,hu^2+gh^2/2,huv]^T$ and $\bm{f}_2 = [hv,huv,hv^2+gh^2/2]^T$. 

We can derive the 1D shallow water equations from the 2D shallow water equations by assuming a rectangular domain with length $L_x$ and $L_y$ in the $x$ and $y$ directions, respectively. If $L_y \ll L_x$ and wall boundary conditions are imposed, then we expect $v$ to be small and $h, u$ to be near-constant along the $y$-direction. These simplifications result in the one-dimensional shallow water equations
\begin{align}
\frac{\partial }{\partial t}\begin{bmatrix}
h\\
hu
\end{bmatrix}+
\frac{\partial }{\partial x}\begin{bmatrix}
hu\\
hu^2+gh^2/2\\
\end{bmatrix}= 0.
\label{eq:SWE_1D}
\end{align}

The mathematical entropy for the shallow water equations corresponds to total energy, and is given by
\[
S(\bm{u}) = \frac{1}{2} \LRp{h \nor{U}^2 + gh^2}
\]
where $\nor{U}^2 = u^2$ in one dimension and  $\nor{U}^2 = u^2 + v^2$ in two dimensions. The entropy variables $\bm{v}$ in two dimensions are given by 
\[
v_1 = gh - \frac{1}{2}\nor{U}^2, \qquad v_2 = u, \qquad v_3 = v.
\]
In one dimension, the entropy variables are simply $\bm{v} = \LRs{v_1,v_2}^T$. The inverse mapping in 2D is given by
\[
h = \frac{v_1 + \frac{1}{2}\nor{U}^2}{g}, \qquad hu = \frac{v_1 + \frac{1}{2}\nor{U}^2}{g} v_2 = h v_2, \qquad hv = \frac{v_1 + \frac{1}{2}\nor{U}^2}{g} v_3 = hv_3.
\]
where we can compute $\nor{U}^2 = v_2^2 + v_3^2$ in terms of the entropy variables. The inverse mapping in 1D follows by ignoring $hv$ and setting $\nor{U}^2 = v_2^2$.

In this paper, we only consider systems where the conservative variables, the entropy potential, and the numerical fluxes can all be transformed between 1D and 2D.
\end{subequations}

%% file: section3.tex
\section{Entropy stable DG discretizations in 1D and 2D}
\label{sec:3}
\begin{subequations}
\numberwithin{equation}{section}
\subsection{On notation}

We adopt a notation which distinguishes between discretized and continuous quantities. 
Unless otherwise specified, continuous vector and matrix quantities are denoted using lower and upper case bold font, respectively. We denote spatially discrete quantities using a bold sans serif font. Finally, the output of continuous functions evaluated over discrete vectors is interpreted as a discrete vector. 

For example, if $\fnt{x}$ denotes a vector of point locations, i.e., $(\fnt{x})_i = \bm{x}_i$, then $u(\fnt{x})$ is interpreted as the vector 
\[	({u}(\fnt{x}))_i = {u}(\bm{x}_i).
\]
Similarly, if $\fnt{u} = {u}(\fnt{x})$, then ${f}(\fnt{u})$ corresponds to the vector
\[
	({f}(\fnt{u}))_i = {f}(u(\bm{x}_i)).
\]
Vector-valued functions are treated similarly. For example, given a vector-valued function $\bm{f}:\mathbb{R}^n\rightarrow \mathbb{R}^n$ and a vector of coordinates $\fnt{x}$, $\LRp{\bm{f}(\fnt{x})}_i = \bm{f}(\bm{x}_i)$.

\subsection{Discretization matrices for high order DG methods}

We construct entropy stable numerical schemes for networks based on high order entropy stable DG formulations in 1D and 2D. These formulations ensure entropy stability over individual segments and subdomains of a multi-dimensional network \cite{chan2018discretely}. We begin by introducing some mathematical notation. We denote the reference element by $\widehat{D}$ with boundary $\partial{\widehat{D}}$. In 1D, the reference element is the interval $[-1,1]$ and in 2D, the reference element is the bi-unit right triangle. We construct entropy conservative schemes on multiple elements, where the domain $\Omega$ is broken up into $K$ non-overlapping elements $D^k$. Each element can be represented as the affine mapping $\Phi_k$ of the reference interval $\widehat{D}$. Because this mapping is affine, $J^k$ (the determinant of the Jacobian of $\Phi_k$) is constant over each element. We use $\widehat{n}_i$ to represent the $i$th component of the outward normal vector scaled by the face Jacobian on the boundary of the reference element. 

The solution is approximated over the reference element by polynomials of total degree $N$
\begin{align}
P^N(\widehat{D}) = \{ \widehat{x}^i\widehat{y}^j,\quad (\widehat{x},\widehat{y}) \in \widehat{D}, \quad 0 \leq i+j \leq N \}.
\end{align}
We denote the dimension of $P^N$ as $N_p = \rm{dim}(P^N(\widehat{D}))$. Moreover, let $\{\phi_i\}_{i=1}^{N_p}$ denote a basis for $P^N$, such that for  $u(\bm{x}) \in P^N(\widehat{D})$,  there exist coefficients $u_i$ such that
\begin{align}
u(\bm{x}) = \sum_{i=1}^{N_p} u_i \phi_i( \widehat{\bm{x}} ), \hspace{1cm} P^N(\widehat{D}) = {\rm span} \{\phi_i(\widehat{\bm{x}}) \}_{i=1}^{N_p}.
\end{align}

We also assume the use of volume and surface quadrature rules $\LRc{\widehat{\bm{x}}_i, w_i}_{i=1}^{N_q}$, $\LRc{\widehat{\bm{x}}^f_i, w^f_i}_{i=1}^{N^f_q}$. We denote the number of volume and surface quadrature nodes by $N_q$ and $N^f_q$ respectively, and assume that the volume and surface rules are exact for polynomials of degree $2N-1$ and $2N$, respectively. We furthermore assume the volume quadrature is sufficiently accurate such that the mass matrix is positive-definite. 

We now introduce quadrature-based operators. Let $\fnt{W}$ and  $\fnt{W}_{f}$ denote diagonal matrices whose entries are volume and surface quadrature weights. We then define the volume and surface quadrature interpolation matrices $\fnt{V}_q$ and $\fnt{V}_f$ as:
\begin{align}
(\fnt{V}_q)_{ij} &= \phi_j(\widehat{\bm{x}}_i),\quad 1 \leq j \leq N_p,\quad  1 \leq i \leq N_q, \\
(\fnt{V}_f)_{ij} &= \phi_j(\widehat{\bm{x}}_i^f),\quad 1 \leq j \leq N_p,\quad  1 \leq i \leq N_q^f,
\end{align}
The matrix $\fnt{V}_q$ maps coefficients of $\fnt{u} = \LRs{u_1, u_2, \ldots, u_{N_p}}$ in terms of polynomial basis to evaluations of $u(\bm{x})$ at volume quadrature points and, similarly, the matrix $\fnt{V}_f$ interpolates $\fnt{u}$ to surface quadrature points. 

We now define $\fnt{D}_i$ as the differentiation matrix with respect to the $i$th coordinate. $\fnt{D}_i$ is defined implicitly by:
\begin{align}
u(\bm{x}) = \sum_{i=1}^{N_p} u_i \phi_i( \widehat{\bm{x}} ), \hspace{1cm} \frac{\partial u}{\partial \widehat{x}_i} = \sum_{j=1}^{N_p} (\fnt{D}_i\fnt{u})_j \phi_j( \widehat{\bm{x}} ).
\end{align}
Here, $\fnt{D}_i$ maps the basis coefficients of some polynomial $u \in P^N$ to coefficients of its $i$th directional derivative with respect to the reference coordinate $\widehat{\bm{x}}_i$.

With the matrix $\fnt{V}_q$, we can now introduce the element mass matrix whose entries are the evaluations of inner products of different basis functions with quadrature points:
\begin{align}
\fnt{M} = \fnt{V}_q^T\fnt{W}\fnt{V}_q, \hspace{0.5cm} \fnt{M}_{ij} = \sum_{k=1}^{N_q} w_k\phi_j(\bm{\widehat{x}}_k)\phi_i(\bm{\widehat{x}}_k) \approx \int_{\widehat{D}} \phi_j\phi_i d\bm{\widehat{x}} = (\phi_j,\phi_i)_{\widehat{D}}.
\end{align}

We can define the quadrature-based $L^2$ projection matrix $\fnt{P}_q$, by inverting the mass matrix:
\begin{align}
\fnt{P}_q = \fnt{M}^{-1}\fnt{V}_q^T\fnt{W}.
\end{align}
The matrix $\fnt{P}_q$ maps a function in terms of its evaluations at quadrature points to its coefficients of the $L^2$ projection in the basis $\phi_i(\bm{\widehat{x}})$.

In d-dimensions, define the following matrices
\begin{align}
\fnt{\widehat{Q}}_i = \fnt{M}\fnt{D}_i, \hspace{1cm} \fnt{B}_i = \fnt{W}_f \diag{\bm{\widehat{n}}_i}, \hspace{1cm} i=1,...,d.
\end{align}
With the above definitions, we have can show that \cite{chan2018discretely}
\begin{align}
\label{eq:SBP_Q}
\fnt{\widehat{Q}}_i + \fnt{\widehat{Q}}_i^T = \fnt{V}_f^T\fnt{B}_i\fnt{V}_f.
\end{align}

By combining the projection matrix $\fnt{P}_q$ with the matrix $\fnt{\widehat{Q}}_i$, we can construct a nodal differentiation operator at quadrature points \cite{chan2018discretely}:
\begin{align}
\fnt{Q}_i = \fnt{P}_q^T\fnt{\widehat{Q}}_i\fnt{P}_q
\end{align}

We also define the the matrix $\fnt{E}$, which extrapolates volume quadrature nodes to surface quadrature nodes, as 
\begin{align}
\fnt{E} = \fnt{V}_f\fnt{P}_q
\end{align}
Then we have the following generalized summation-by-parts (SBP) property:
\begin{align}
\fnt{Q}_i + \fnt{Q}_i^T = \fnt{E}^T\fnt{B}_i\fnt{E}  
\end{align}

Entropy stable formulations for nonlinear conservation laws usually introduce numerical flux terms which couple together all degrees of freedom on neighboring elements \cite{crean2017high}.  To avoid this, we introduce the hybridized operator $\fnt{Q}_{i,h}$, which is given explicitly as 
\begin{align}
\fnt{Q}_{i,h} = \frac{1}{2}
\begin{bmatrix}
\fnt{Q}_i - \fnt{Q}_i^T & \fnt{E}^T\fnt{B}_i\\
-\fnt{B}_i\fnt{E} & \fnt{B}_i\\
\end{bmatrix}.
\end{align}
This operator is designed to be applied to vectors of solution values at both volume and surface quadrature nodes and mimics the structure of boundary terms used in hybridized DG methods \cite{chen2019review}.  When used in a DG formulation, it allows one to construct entropy stable formulations using more standard DG numerical fluxes. We have the following theorem:
\begin{Thm}
$\fnt{Q}_{i,h}$ satisfies the $SBP-like$ property \cite{chan2018discretely}:
\begin{align}
\fnt{Q}_{i,h} +\fnt{Q}_{i,h}^T = \fnt{B}_{i,h}, \hspace{1cm} \fnt{B}_{i,h} = \begin{bmatrix}
\fnt{0} &  \\
  &  \fnt{B}_i\\
\end{bmatrix},
\end{align}
 and $\fnt{Q}_{i,h} \fnt{1} = 0$, where $\fnt{1}$ is the vector of all ones
\label{thm:SBP_Q_h^i}
\end{Thm}
We also construct differentiation and boundary matrices $\fnt{Q}^k_{i,h},\fnt{B}^k_{i}$ on the physical element $D^k$ through the chain rule and linear combinations of differentiation matrices on the reference element. These will be used to construct entropy conservative and entropy stable schemes in the following section. Let $g^k_{ij} = J^k\pd{x_i}{\hat{x}_j}$ denote geometric terms on $D^k$. Then, physical SBP operators can be constructed by taking a linear combination of the reference SBP operators
\begin{equation}
\fnt{Q}^k_{i,h} = \sum_{j=1}^d g^k_{ij} \fnt{Q}_{i,h}.
\label{eq:physSBP}
\end{equation}
It was shown in \cite{chan2018discretely} that these operators satisfy a physical SBP property
\[
\fnt{Q}^k_{i,h} + \LRp{\fnt{Q}^k_{i,h}}^T = \begin{bmatrix}
\fnt{0} &  \\
  &  \fnt{B}^k_{i}\\
\end{bmatrix},
\]
where $\fnt{B}^k_{i}$ is a diagonal matrix containing the $i$th component of the outward unit normal on $D^k$ scaled by quadrature weights and surface Jacobian factors. We note that this construction uses the fact that $g^k_{ij}$ is constant over $D^k$ for affinely mapped elements. We also introduce the physical mass matrix $\fnt{M}^k = J^k \fnt{M}$, which is scaled by a Jacobian factor. Finally, we note that physical SBP operators on curved elements can be constructed using a ``split form'' version of (\ref{eq:physSBP}) \cite{chan2019discretely}.

Finally, we introduce $\bm{v}_h$ as the $L^2$ projection of the entropy variables and $\bm{\tilde{u}}$ as the evaluations of the conservative variables in terms of the $L^2$ projected entropy variables
\begin{align}
\bm{u}_q = \fnt{V}_q\bm{u}_h, \hspace{1cm} \bm{v}_q = \bm{v}(\bm{u}_q), \hspace{1cm} \bm{v}_h = \fnt{P}_q\bm{v}_q, 
\end{align}
\begin{align}
\bm{\tilde{v}} = 
\begin{bmatrix}
\bm{\tilde{v}}_q\\ \bm{\tilde{v}}_f
\end{bmatrix}
 = \begin{bmatrix}
 \fnt{V}_q\\ \fnt{V}_f
 \end{bmatrix}
 \bm{v}_h, \hspace{1cm}
 \bm{\tilde{u}} = 
\begin{bmatrix}
\bm{\tilde{u}}_q\\ \bm{\tilde{u}}_f
\end{bmatrix} = 
 \bm{u}(\bm{\tilde{v}}).
\end{align}
Here $\bm{u}_q$ and $\bm{v}_q$ denote the conservative variables and entropy variables evaluated at the volume quadrature points. The vector $\bm{\tilde{v}}$ denotes the evaluations of the $L^2$ projection of the entropy variables at both volume and surface quadrature points, while $\bm{\tilde{u}}$ denotes the evaluations of the conservative variables in terms of the projected entropy variables $\bm{u}(\Pi_N \bm{v})$, where $\Pi_N$ denotes the $L^2$ projection operator. 

\subsection{Entropy conservation and  flux differencing}\label{sec:EC_flux}
In this section, we introduce entropy conservative numerical fluxes and discrete formulations \cite{fisher2013high, gassner2018br1, chen2017entropy, gassner2016split}. To construct entropy stable schemes in $d$ dimensions, we utilize entropy conservative fluxes as defined in \cite{tadmor1987numerical}:
\begin{Def}{}
\label{def:consevative_flux}
Let $\bm{f}_{i,S}(\bm{u}_L, \bm{u}_R)$ be a bivariate function which is symmetric and consistent with the flux function $\bm{f}_i(\bm{u})$, for $i = 1,...,d$
\begin{align}
\bm{f}_{i,S}(\bm{u},\bm{u}) = \bm{f}_i(\bm{u}), \hspace{1cm} \bm{f}_{i,S}(\bm{u},\bm{v}) = \bm{f}_{i,S}(\bm{v},\bm{u}).
\label{eq:ec_def_1}
\end{align}
$\bm{f}_{i,S}(\bm{u}_L, \bm{u}_R)$ is called entropy conservative if, for entropy variables $\bm{v}_L = \bm{v}(\bm{u}_L)$ and $\bm{v}_R = \bm{v}(\bm{u}_R)$, 
\begin{align}
\left( \bm{v}_L - \bm{v}_R\right)^T \bm{f}_{i,S}\left(\bm{u}_L,\bm{u}_R\right) = \psi_i(\bm{v}(\bm{u}_L)) - \psi_i(\bm{v}(\bm{u}_R)).
\label{eq:ec_def_2}
\end{align}
\end{Def}
The flux $\bm{f}_{i,S}$ can be used to construct entropy conservative and entropy stable finite volume methods. 
This numerical flux can also be used to construct discretely entropy stable DG schemes using an approach referred to as flux differencing \cite{fisher2013high, carpenter2014entropy,gassner2018br1, chen2017entropy}. 

Using flux differencing \cite{chan2018discretely, crean2018entropy}, we can approximate the derivative of $\bm{f}_i(\bm{u}(x))$ using the differentiation matrices $\fnt{Q}^k_{i,h}$ and $\bm{f}_{i,S}$. We define a flux matrix $\fnt{F}_i$ by evaluating $\bm{f}_{i,S}$ using solution values at quadrature points:
\begin{align}
(\fnt{F}_i)_{lm} = \bm{f}_{i,S}(\bm{u}_l, \bm{u}_m), \qquad 1 \leq l,m \leq N_q.
\label{eq:flux_matrix}
\end{align}
Then, $\int \frac{\partial \bm{f}_i(\bm{u})}{\partial x_i}$ can be approximated by the term $2(\fnt{Q}^k_{i,h}\circ\fnt{F}_i)\bm{1}$, where $\circ$ denotes the Hadamard product between two matrices. A discrete entropy conservative formulation is then given as follows on an element $D^k$:
\begin{align}
\fnt{M}^k\td{\fnt{u}}{t} + \sum_{i=1}^d \begin{bmatrix}
\fnt{V}_q\\\fnt{V}_f
\end{bmatrix}^T\LRp{2\fnt{Q}^k_{i,h}\circ \fnt{F}_i^k}\fnt{1} + \fnt{V}_f^T\fnt{B}^k_{i}\LRp{\bm{f}_{i,S}(\tilde{\fnt{u}}^{k+},\tilde{\fnt{u}}^k)-\bm{f}_i(\tilde{\fnt{u}}_f^k)} = \fnt{S},\label{eq:dgform}
\end{align}
where $\fnt{S}$ denotes source terms depending on the bottom geometries. 

Entropy stable schemes can be constructed from an entropy conservative scheme by adding appropriate penalization term that dissipate entropy at element interfaces \cite{chen2017entropy}. This modification converts schemes which satisfy a global entropy equality into schemes which satisfy a global entropy inequality. 
\end{subequations}

%% file: section4.tex
\section{Entropy stable coupling terms for junctions and multi-dimensional interfaces}
\label{sec:4}
\begin{subequations}
\numberwithin{equation}{section}



In this section, we describe how to construct entropy conservative and entropy stable coupling terms for 1D junctions (which we will refer to as 1D-1D couplings) and couplings between 1D and 2D domains (which we will refer to as 1D-2D couplings). We refer to both individual channels in networks and 1D and 2D parts of a multi-dimensional domain as subdomains, and discretize over each individual sub-domain using the entropy stable DG methods in \cite{chan2018discretely}. We note, however, that the proposed approach is applicable to any entropy stable SBP-type discretization. 

We proceed by designing interface coupling terms which are entropy conservative in the sense that they do not contribute to entropy production. The addition of entropy dissipation (both at coupling interfaces and over each sub-domain) then yields entropy stable schemes over networks or multi-dimensional domains. 
\subsection{Notation and assumptions}
We first introduce notation for 1D and 2D fluxes. Let $\bm{f}_{i,S}$ denote the 2D flux, where the subscript $i$ denotes the $i$th component of the 2D fluxes, (where $i=1,2$ correspond to the $x$ and $y$ direction respectively). We denote the 1D flux vector by $\bm{f}_S$, without the subscript for the coordinate index. We use $\bm{f}_{J,1D}$ and $\bm{f}_{J,2D}$ to denote the numerical flux at a multi-dimensional junction on 1D and 2D sides respectively. We use $\bm{f}_{J,i}$ to denote the fluxes at a 1D-1D junction, where the subscript $i$ denotes the numerical flux for the $i$th 1D channel. We use $n_{1D}$ and $n_{2D}$ to denote the sign of the outward surface normal on the 1D and 2D domains at a 1D-2D junction. 

Finally, we restrict ourselves to velocity-based systems, where in $d$ dimensions, the conservative variables contain a velocity or momentum vector with $d$ components. We also assume that 1D system can be derived from 2D system by ignoring or projecting certain velocity components. We note our approach is not directly applicable to systems which do not satisfy these assumptions, such as Burger's equation. 
\subsection{1D-2D domain coupling}
\label{subsec:1D-2D}
\begin{figure}
\begin{center}
\includegraphics[width=.5\textwidth]{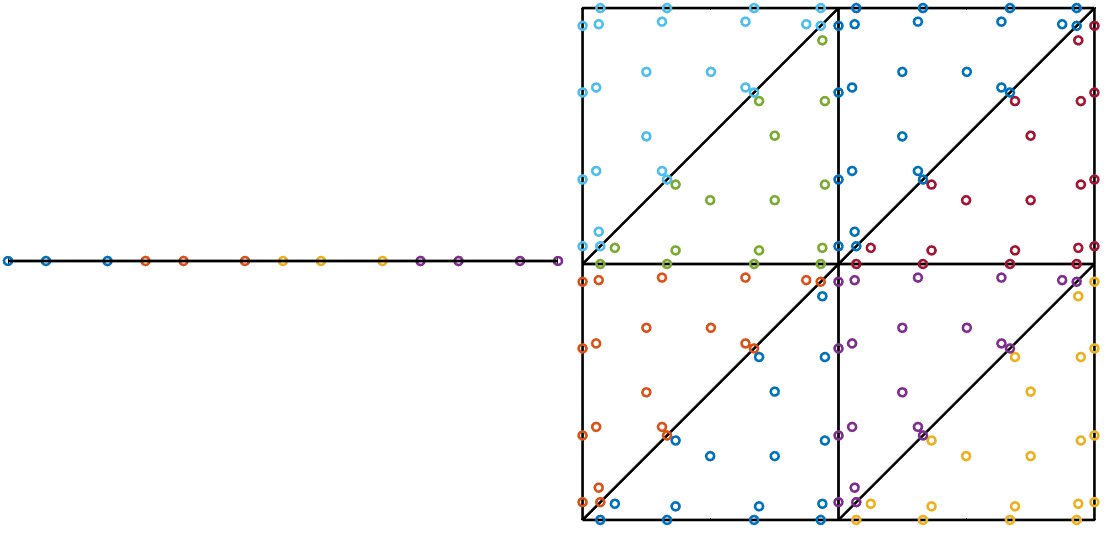}
\caption{A 1D-2D mesh with surface and volume quadrature points overlaid.}
\label{fig:1_D_2D_mesh}
\end{center}
\end{figure}
For flows over networks of 1D channels, one approach to model junction behavior is to explicitly discretize the junction using a 2D mesh \cite{neupane2015discontinuous, bellamoli2018numerical}. Since each junction can be accurately modeled using a relatively small number of elements, this approach inherits most of the computational advantages of a purely 1D network simulation. 
We illustrate the construction of a 1D-2D coupling scheme using the simplified setting shown in Figure~\ref{fig:1_D_2D_mesh}. 
We assume that the 2D domain has an interface which is perpendicular to the 1D domain. This property can be enforced when constructing meshes for 2D domains. We couple 1D and 2D domains together by treating the 1D domain as a channel with the same width of the 2D domain. 

In both 1D and 2D, each cell communicates through fluxes (calculated in terms of ``interior'' and ``exterior'' values of the conservative variables) at the interfaces, and we treat the coupling interface similarly. We first compute the fluxes at the interface of the coupling in the 2D domain with the following steps. Let $\bm{u}_{2D}$ denote the vector of the conservative variables in the ``interior'' of 2D cell, and let $\bm{u}^+_{2D}$ denotes the ``exterior'' values of the conservative variables from the neighbouring cell. 
The conservative variables in 1D have fewer fields than the variables in 2D. To resolve this mismatch, we introduce a matrix $\bm{R}$ to transform between 1D and 2D variables. For shallow water equations in 1D and 2D, we have
\begin{align}
\bm{R} = \begin{bmatrix}
1 &0\\
0 &n_1\\
0 &n_2\\
\end{bmatrix}, \qquad \bm{R}^T\bm{R} = \bm{I},
\label{eq:1d2d_transformation_matrix}
\end{align}
where $[n_1,n_2]^T$ is the unit outward normal vector for the 1D domain with respect to the 2D coordinates. Then we can use $\bm{R}$ to transform between 1D and 2D variables as follows:

\begin{align}
\bm{R}\bm{u}_{1D}&= \begin{bmatrix}
1 &0\\
0 &n_1\\
0 &n_2\\
\end{bmatrix}
\begin{bmatrix}
h\\
hu_{1D}\\
\end{bmatrix} \Longrightarrow 
\begin{bmatrix}
h\\
hu_{2D}\\
hv_{2D}\\
\end{bmatrix}, \\ 
\bm{R}^T\bm{u}_{2D}&= \begin{bmatrix}
1 &0 &0\\
0 &n_1 &n_2\\
\end{bmatrix}
\begin{bmatrix}
h\\
hu_{2D}\\
hv_{2D}\\
\end{bmatrix} \Longrightarrow
\begin{bmatrix}
h\\
hu_{1D}\\
\end{bmatrix}.
\label{eq:1d2d_transform}
\end{align}

Note that $\bm{R}^T\bm{u}_{2D}$ constructs the $1D$ momentum from the normal component of the 2D momentum. Similarly, we can apply this $\bm{R}$ operator to the entropy conservative fluxes and entropy variables in both 1D and 2D. 

Let $\bm{u}_{f,2D}$ and $\bm{u}_{f,1D}$ denote the vector of values of the conservative variables at surface quadrature points in the 2D and 1D domains on general elements and let $\bm{u}_{J,2D}$ and $\bm{u}_{J,1D}$ denote the vector of values of the conservative variables at surface quadrature points on the 2D and 1D side of a 1D-2D interface. Then the flux on the 2D side of the 1D-2D interface $\bm{f}_{J,2D}$ is defined as 
\begin{align}
\bm{f}_{J,2D}(\bm{u}_{J,2D}^+, \bm{u}_{J,2D}) = \sum_{i = 1}^d n_i\bm{f}_{i,S}(\bm{u}_{J,2D}^+, \bm{u}_{J,2D}), \qquad \bm{u}_{J,2D}^+ = \bm{R}\bm{u}_{J,1D}.
\label{eq:1D-2D_flux2D}
\end{align}
We next apply the transformation matrix $\bm{R}^T$ to the 2D fluxes to reduce them down to 1D dimension.
We define the surface flux for 1D side of the 1D-2D interface as follows:
\begin{align}
\bm{f}_{J,1D} = \frac{\fnt{w}_{J,f}^T\LRp{\bm{R}^T\bm{f}_{J,2D}}}{\fnt{w}_{J,f}^T\fnt{1}}n_{1D}, 
\label{eq:1D-2D_flux}
\end{align}
where $\fnt{w}_{J,f}$ is the vector containing the quadrature weights of the surface quadrature points on the entire 1D-2D interface.
Note that operation $\bm{R}^T\bm{f}_{J,2D}$ is performed point-wise at each surface quadrature point.

We now present a proof of entropy conservation for our 1D-2D coupling scheme. Let $K_{1D}$ and $K_{2D}$ denote the number of elements in the 1D and 2D domains respectively. Let $J^{k}_{1D}$ and $J^{k}_{2D}$ be the Jacobian on the $k$th elements in the 1D and 2D domains respectively. We first introduce the following lemma:
\begin{lemma}
Let $\fnt{\tilde{u}}_{1D}$ and $\fnt{\tilde{u}}_{2D}$ denote the 1D and 2D entropy projected conservative variables, and let $\fnt{r}_{1D}(\fnt{\tilde{u}}_{1D})$ and $\fnt{r}_{2D}(\fnt{\tilde{u}}_{2D})$ denote the entropy stable DG spatial formulation such that
\begin{align}
&\fnt{r}_{1D}(\fnt{\tilde{u}}_{1D}) = \sum_{k=1}^{K_{1D}} 
(\begin{bmatrix}
\fnt{V}_q\\\fnt{V}_f
\end{bmatrix}^T
2\LRp{\LRp{\fnt{Q}^k - (\fnt{Q}^k)^T}\circ \fnt{F}}\fnt{1} + \fnt{V}_f^T\fnt{B}\bm{f}_{1D}^*),\\
&\fnt{r}_{2D}(\fnt{\tilde{u}}_{2D}) = \sum_{k=1}^{K_{2D}} \sum_{i=1,2} (\begin{bmatrix}
\fnt{V}_q\\\fnt{V}_f
\end{bmatrix}^T
2\LRp{\fnt{Q}^k_i - (\fnt{Q}^k_i)^T\circ \fnt{F}_i}\fnt{1} + \fnt{V}_f^T\fnt{B}^k_i\bm{f}_{i,2D}^*),\\
&(\fnt{F})_{j,k} = \bm{f}_S(\fnt{\tilde{u}}_j, \fnt{\tilde{u}}_k), \qquad (\fnt{F}_i)_{j,k} = \bm{f}_S(\fnt{\tilde{u}}_j, \fnt{\tilde{u}}_k),
\end{align}
where 
\begin{align}
\bm{f}_{1D}^* &= \bm{f}_S(\fnt{\tilde{u}}_{f,1D}^+, \fnt{\tilde{u}}_{f,1D})\\
\bm{f}_{i,2D}^* &= \bm{f}_{i,S}(\fnt{\tilde{u}}_{f,2D}^+, \fnt{\tilde{u}}_{f,2D})\\
\end{align}
on interior interfaces and 
\begin{align}
\bm{f}_{1D}^* = 0 \qquad \bm{f}_{i,2D}^* =0
\end{align}
at the junction boundaries of the 1D and 2D domains. Let $\fnt{v}_{1D}$ and $\fnt{v}_{2D}$ denote the projected entropy variables. Assuming that the entropy flux $\psi(\tilde{\fnt{u}}) - \fnt{v}^T \bm{f}^* = 0$ on non-junction boundaries, then
\begin{align}
\fnt{v}_{1D}^T\fnt{r}_{1D}(\fnt{\tilde{u}}_{1D}) &= \psi(\fnt{\tilde{u}}_{J,1D})\\
\fnt{v}_{2D}^T\fnt{r}_{2D}(\fnt{\tilde{u}}_{2D}) &= \sum_{i=1,2}\psi_i(\fnt{\tilde{u}}_{J,2D})n_i.
\end{align}
\label{lemma:entropy_protential1D2D}
\end{lemma}

The proof of the lemma is restatement of results from \cite{chan2018discretely}. Next, assume a simple 1D-2D coupling as shown in Figure \ref{fig:1_D_2D_mesh}, where we assume the 1D and 2D domain are channels of the same width denoted by $A$. The 1D-2D junction terms can be expressed as
\begin{align}
\sum_{k=1}^{K_{1D}} \fnt{M}^k \td{\fnt{u}_{1D}}{t} + \fnt{r}_{1D}(\fnt{\tilde{u}}_{1D}) + \fnt{V}_f^T \bm{f}_{J,1D}(\fnt{\tilde{u}}_{J,1D}, \fnt{\tilde{u}}_{J,2D})n_{1D} = 0\\
\sum_{k=1}^{K_{2D}}  \fnt{M}^k \td{\fnt{u}_{2D}}{t} + \fnt{r}_{1D}(\fnt{\tilde{u}}_{2D}) + \fnt{V}_f^T\fnt{W}_{f} \bm{f}_{J,2D}(\fnt{\tilde{u}}_{J,1D}, \fnt{\tilde{u}}_{J,2D})n_{2D} = 0,
\end{align}
where $\bm{f}_{J,2D}(\fnt{\tilde{u}}_{J,1D}, \fnt{\tilde{u}}_{J,2D})$ and 
$\bm{f}_{J,1D}(\fnt{\tilde{u}}_{J,1D}, \fnt{\tilde{u}}_{J,2D})$ are defined in (\ref{eq:1D-2D_flux}) and (\ref{eq:1D-2D_flux2D}) respectively. Together, these 1D-2D discretizations imply entropy conservation: 
\begin{theorem}
Let fluxes at a 1D-2D interface be defined as in (\ref{eq:1D-2D_flux2D}) and (\ref{eq:1D-2D_flux}). Let $\fnt{u}_{q,1D}$ and $\fnt{u}_{q,2D}$ denote the values of the 1D and 2D solutions at quadrature points on the $k$th 1D or 2D element. Then, assuming continuity in time and that the entropy flux $\psi(\tilde{\fnt{u}}) - \fnt{v}^T \bm{f}^* = 0$ on  non-junction boundaries, the DG scheme defined by (\ref{eq:dgform}) is entropy conservative in the sense that
\begin{align}
A\sum_{k=1}^{K_{1D}}\bm{1}^TJ^{k}_{1D}\fnt{W}\frac{{\rm{d}}S(\fnt{u}_{q,1D})}{{\rm{d}}t} + \sum_{k=1}^{K_{2D}}\bm{1}^TJ^{k}_{2D}\fnt{W}\frac{{\rm{d}}S(\fnt{u}_{q,2D})}{{\rm{d}}t} =0,
\end{align}
which is a quadrature approximation to 
\begin{align}
\frac{\partial}{\partial t}\LRp{A\int_{\Omega_{1D}}S(\bm{u})+\int_{\Omega_{2D}}S(\bm{u})}=0,
\end{align}
where $\Omega_{1D}$ and $\Omega_{2D}$ denote the 1D and 2D domain respectively. 
\label{thm:entropy_conservative}
\end{theorem}
\begin{proof}
It is sufficient to prove entropy conservation for the setup shown in Figure \ref{fig:1_D_2D_mesh}. From the results in \cite{chan2018discretely} and Lemma \ref{lemma:entropy_protential1D2D}, we only need to show that the flux contributions from 1D and 2D sides of the 1D-2D junction interface cancel. Testing with the entropy variables, scaling with the Jacobian $J^k$ in both domains and width $A$ of the 1D domain, it can be shown that the 1D and 2D schemes each satisfy
\begin{align}
A\sum_{k=1}^{K_{1D}}\bm{1}^TJ^{k}_{1D}\fnt{W}\frac{{\rm{d}}S(\fnt{u}_{q,1D})}{{\rm{d}}t} &=A\bm{1}^Tn_{1D}\LRp{ \psi_
{1D}(\fnt{\tilde{u}}_{J,1D})-\fnt{\tilde{v}}_{1D}^T\bm{f}_S(\fnt{\tilde{u}}_{J,1D}^+, \fnt{\tilde{u}}_{J,1D})},
\label{eq:1d2d_1d_ec}
\end{align}
\begin{align}
\sum_{k=1}^{K_{2D}}\bm{1}^TJ^{k}_{2D}\fnt{W}\frac{{\rm{d}}S(\fnt{u}_{q,2D})}{{\rm{d}}t} &= \sum_{i=x,y} \bm{1}^T\fnt{w}_fn_{i,2D}\LRp{ \psi_i(\fnt{\tilde{u}}_{J,2D})-\fnt{\tilde{v}}_{2D}^T\bm{f}_{i,S}(\fnt{\tilde{u}}_{J,2D}^+, \fnt{\tilde{u}}_{J,2D})}.
\label{eq:1d2d_2d_ec}
\end{align}
Without loss of generality, we assume $n_{1,2D} = 1$ and $n_{2,2D} = 0$, so we denote $n_{1,2D}$ with shorthand $n_{2D}$. On 1D side of the junction, the flux contribution is:
\begin{align}
&A n_{1D}\fnt{(\tilde{v}}_{J,1D}^+)^T\bm{f}_{J,1D} = A ^Tn_{1D}\fnt{(\tilde{v}}_{J,1D}^+)^T\frac{\fnt{w}_{J,f}^T\LRp{\bm{R}^T\bm{f}_{J,2D}}}{\fnt{w}_{J,f}^T\bm{1}},
\end{align}
where $n_{1D}$ is the normal in the 1D domain. Using that $\fnt{w}_{J,f}^T\fnt{1} = A$, we have the flux contribution from the 1D side is 
\begin{align}
&A n_{1D}(\fnt{\tilde{v}}_{J,1D}^+)^T\frac{\fnt{w}_f^T\LRp{\bm{R}^T\bm{f}_{J,2D}}}{\fnt{w}_{J,f}^T\bm{1}}\\
= &n_{1D}(\fnt{\tilde{v}}_{J,1D}^+)^T\fnt{w}_{J,f}^T\bm{R}^T\bm{f}_{J,2D}\LRp{(\bm{R}{\fnt{\tilde{u}}_{J,1D}})\bm{1}, \fnt{\tilde{u}}_{J,2D})}\label{eq:1d2d-1dterms_prev}\\
= &n_{1D}(\bm{R}\fnt{\tilde{v}}_{J,1D}^+\bm{1})^T\fnt{W}_{J,f}\bm{f}_{J,2D}\LRp{(\bm{R}{\fnt{\tilde{u}}_{J,1D}})\bm{1}, \fnt{\tilde{u}}_{J,2D})}
\label{eq:1d2d-1dterms}
\end{align}
where, $\fnt{W}_{J,f}$ is a diagonal matrix that contains the surface quadrature weights at the junction interface. 
In going from (\ref{eq:1d2d-1dterms_prev}) to (\ref{eq:1d2d-1dterms}), we use that $\fnt{w}_{J,f}^T = \bm{1}^T\fnt{W}_{J,f}$ by the definitions of $\fnt{W}_{J,f},\fnt{w}_{J,f}$. We also use that since multiplication by $\bm{R}^T$ acts on discrete solution values on junction surface quadrature points (to account for the fact that $\fnt{\tilde{v}}_{J,1D}^+$ is a scalar and $\bm{f}_{J,2D}$ is a vector), multiplication by $\fnt{W}_f^T$ and $\bm{R}^T$ commute.

Similar in 2D, we have the following flux contribution on the junction
\begin{align}
n_{2D}(\fnt{\tilde{v}}_{J,2D}^+)^T\fnt{W}_{J,f}\bm{f}_{J,2D}\LRp{(\bm{R}{\fnt{\tilde{u}}_{J,1D}})\bm{1}, \fnt{\tilde{u}}_{J,2D})}
\label{eq:1d2d-2dterms}
\end{align}
Because the outward normal for the 1D domain is the negative of the outward normal for the 2D domain, we can combine the 1D and 2D surface terms from (\ref{eq:1d2d-1dterms}) and (\ref{eq:1d2d-2dterms}) together. Using that $\fnt{1}^T\fnt{w}_{J,f} = A$, we have 
\begin{align}
&n_{1D}(\bm{R}\fnt{\tilde{v}}_{J,1D}^+\bm{1})^T \fnt{W}_{J,f} \bm{f}_{J,2D}+ n_{2D}(\fnt{\tilde{v}}_{J,2D}^+)^T\fnt{W}_{J,f}\bm{f}_{J,2D}\\
=& \LRp{(\bm{R}\fnt{\tilde{v}}_{J,1D}^+\bm{1})^T - (\fnt{\tilde{v}}_{J,2D}^+)^T}\fnt{W}_{J,f} \bm{f}_{J,2D}\LRp{(\bm{R}{\fnt{\tilde{u}}_{J,1D}})\bm{1}, \fnt{\tilde{u}}_{J,2D})}
\end{align}
First, note $\fnt{\tilde{u}}_{J,2D} = \fnt{u}(\fnt{\tilde{v}}_{J,2D})$ by construction. Then, note that the 2D mapping between conservative and entropy variables reduces to a 1D mapping in the normal direction upon multiplication by $\bm{R}^T$. Thus,  $\fnt{u}\LRp{(\bm{R}\fnt{v}_{1D}}\bm{1})=\bm{R}\fnt{\tilde{u}}_{1D}\bm{1}$,  where $\fnt{u}_{1D} = \fnt{u}(\fnt{v}_{1D})$. 

The flux contributions from 1D and 2D side of the junction cancel since the operator $R$ transform the 1D projected entropy variables to 2D, such that $(\bm{R}\fnt{\tilde{v}}_{J,1D}^+\bm{1})^T = (\fnt{\tilde{v}}_{J,2D}^+)^T$. Therefore, with $\psi(\tilde{\fnt{u}}) - \fnt{v}^T \bm{f}^* = 0$ on non-junction boundaries, we have
\begin{align}
A\sum_{k=1}^{K_{1D}}\bm{1}^TJ^{k}\fnt{W}\frac{{\rm{d}}S(\fnt{u}_{q,1D})}{{\rm{d}}t} + \sum_{k=1}^{K_{2D}}\bm{1}^TJ^{k}\fnt{W}\frac{{\rm{d}}S(\fnt{u}_{q,2D})}{{\rm{d}}t} =0.
\end{align}
\end{proof}

We note that a 2D domain can be coupled to multiple 1D domains. For example, in Figure \ref{fig:1_D_2D_mesh_split}, we calculate the flux for Channel 1 with fluxes from the top two cells and the flux for Channel 2 with the bottom six cells of the 1D-2D interface.
\begin{figure}[H]
\begin{center}
\includegraphics[width=.75\textwidth]{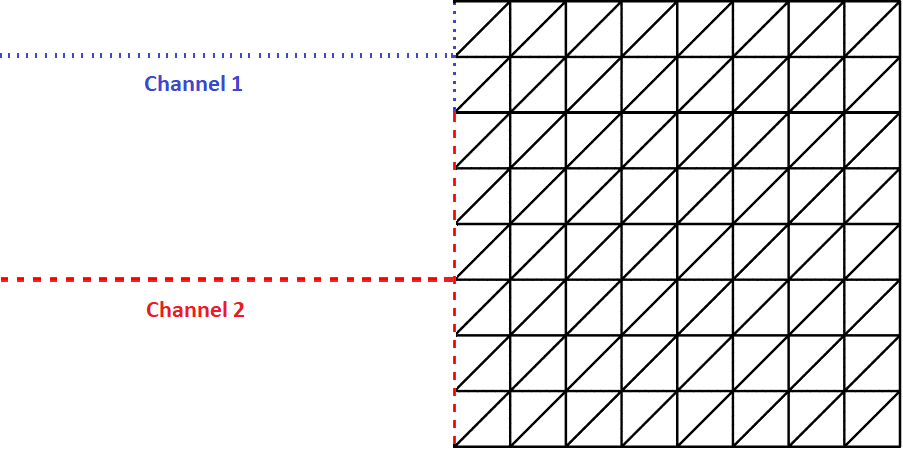}
\caption{Example of a 1D-2D mesh with multiple channels. Channel 1 with its connecting interface are shown using blue dots and channel 2 with its connecting interface are shown using red dashes.}
\label{fig:1_D_2D_mesh_split}
\end{center}
\end{figure}

\subsection{1D-1D junction coupling}
In this section, we introduce entropy stable methods to join multiple (2 or more) 1D domains together without explicitly representing the junction as a 2D domain. This approach simplifies meshing and further reduces the overall computational cost, though it may not accurately reproduce multi-dimensional effects in certain flow setting \cite{bellamoli2018numerical, neupane2015discontinuous}.
\begin{figure}[H]
\begin{center}
\includegraphics[width=.45\textwidth]{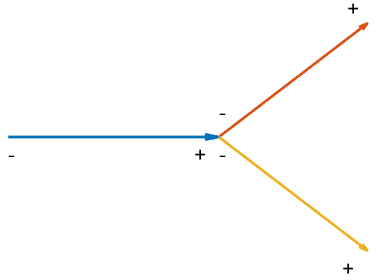}
\caption{1D domains with directions and signs of the outward unit normal for junction.}
\label{fig:1_D_mesh_direction}
\end{center}
\end{figure}
We first present a general method to couple an arbitrary number of 1D domains at a junction. The key is to distribute shared fluxes between each domain connected to the junction in a conservative fashion. Let $I_J$ be a set of all domains at the junction and let $|I_J|$ denotes the dimension of the set. Let $A_i$ denote the width of the $i$th domain at a junction. 
Let $n_{J,i}= \pm{1}$ denote the outward surface normal for channel $i \in I_J$, as shown in Figure \ref{fig:1_D_mesh_direction}. We introduce on channel $i \in I_J$ the values of solution $\bm{u}_{f,i}$ and  junction flux $\bm{f}_{J,i}$

\begin{align}
\bm{f}_{J,i}(\bm{u}_{f,i}^+, \bm{u}_{f,i})= \sum_{j\in I_J} c_{ij} \bm{f}_S(\bm{u}_{f,i}, \bm{u}_{f,j}),
\label{eq:1D_1D_flux}
\end{align}

where we have the weighting coefficients $c_{ij}$. We require these coefficients to satisfy 
\begin{align}
\sum_{j\in I_J} c_{ij} =1, \textrm{ and } A_i c_{ij} = A_j c_{ji},  \textrm{ for all } i,j.
\label{eq:1D_1D_flux_coefficients}
\end{align}
We introduce the convention that non-zero values of $c_{ii}$ represent ``partial" wall boundary conditions (which we describe in more detail in Section \ref{subsec:partial_wall}), where the solution of the a domain partially interacts with itself. 

Let $K_i$ denote the number of elements in the $i$th channel and let $J^{k}_{i}$ be the Jacobian of the $k$th element of the $i$th channel. With Lemma \ref{lemma:entropy_protential1D2D}, we have the following theorem:
\begin{theorem}
Let fluxes at the 1D-1D junction interface be defined as in (\ref{eq:1D_1D_flux}) with the coefficients satisfying (\ref{eq:1D_1D_flux_coefficients}), and let $\fnt{u}_{q,i}$ denote the solution values on the quadrature points on the $i$th channel. Then, assuming continuity in time, the DG scheme defined by (\ref{eq:dgform}) for the entire network is entropy conservative for periodic boundary conditions
\begin{align}
\sum_{i\in I_J} A_i\sum_{k=1}^{K_i}\bm{1}^TJ^{k}_{i}\fnt{W}\frac{{\rm{d}}S(\fnt{u}_{q,i})}{{\rm{d}}t} = 0.
\end{align}
This is a quadrature approximation to 
\begin{align}
\td{}{t}\LRp{\sum_{i\in I_J}A_i\int_{\Omega_i}S(\bm{u}_i)}=0.
\end{align}
\label{thm:entropy_conservative_1D_1D}
\end{theorem}
\begin{proof}
Consider a network of 1D channels which join at a junction. Since we are using an entropy conservative DG formulation on each channel, our scheme is entropy conservative up to the coupling interface, and satisfies on each channel \cite{chan2018discretely}
\begin{align}
\bm{1}^T\fnt{W}\frac{{\rm{d}}S(\fnt{u}_{q,i})}{{\rm{d}}t} = \bm{1}^Tn_{J,i}\LRp{ \psi(\fnt{\tilde{u}}_{f,i})-\fnt{\tilde{v}}_i^T\bm{f}_S},  \qquad \forall i\in I_j.
\end{align}
Scaling by the Jacobian $J^{k}_{i}$ and width, then summing over all channels and elements gives
\begin{align}
\sum_{i\in I_J} A_i\sum_{k=1}^{K_i}\bm{1}^TJ^{k}_{i}\fnt{W}\frac{{\rm{d}}S(\fnt{u}_{q,i})}{{\rm{d}}t} = \sum_{i\in I_J} A_i\sum_{k=1}^{K_i} \bm{1}^Tn_{J,i}\LRp{ \psi(\fnt{\tilde{u}}_{f,i})-\fnt{\tilde{v}}_{J,1D}^T\bm{f}_S(\fnt{\tilde{u}}_{f,i}^+, \fnt{\tilde{u}}_{f,i})}.
\end{align}
Notice that $n^+ = -n$, where $n$ is the outward normal and recall that the entries of $B$ correspond to the value of $n$ at the channel. Then splitting the interface contribution between neighboring elements gives
\begin{align}
\nonumber
-\sum_{i\in I_J} A_i\sum_{k=1}^{K_i} \bm{1}^Tn_{J,i}\LRp{ \fnt{\tilde{v}}_i^T\bm{f}_S(\fnt{\tilde{u}}_{f,i}^+, \fnt{\tilde{u}}_{f,i})} = &-\sum_{i\in I_J} A_i\sum_{k\notin J} \bm{1}^Tn_{J,i}\LRp{ \fnt{\tilde{v}}_i^T\bm{f}_S(\fnt{\tilde{u}}_{f,i}^+, \fnt{\tilde{u}}_{f,i})} \\ 
&- \sum_{i\in I_J} A_i\sum_{k\in J} \bm{1}^Tn_{J,i}\LRp{ \fnt{\tilde{v}}_i^T\bm{f}_{J,i}(\fnt{\tilde{u}}_{f,i}^+, \fnt{\tilde{u}}_{f,i})}.
\label{eq:1d_1d_junction_ec_pf_split}
\end{align}

The first part of (\ref{eq:1d_1d_junction_ec_pf_split}) corresponds to the elements away from junction. By Lemma~\ref{lemma:entropy_protential1D2D}, we only need to consider the numerical flux contribution at the junction. Without loss of generality, assume that only one element from each domain is connected to the junction. For elements at the junction, we can write the remaining term in (\ref{eq:1d_1d_junction_ec_pf_split}) as
\begin{align}
&-\sum_{i\in I_J} A_i\sum_{k\in J} \bm{1}^Tn_{J,i}\LRp{ \fnt{\tilde{v}}_i^T\bm{f}_{J,i}(\fnt{\tilde{u}}_{f,i}^+, \fnt{\tilde{u}}_{f,i})}\\ = &\frac{1}{2}\sum_{i\in I_J} A_i \bm{1}^Tn_{J,i}\LRp{ \sum_{j=1}^{|I_J|}\LRp{\fnt{\tilde{v}}_j-\fnt{\tilde{v}}_i}^Tc_{ij}\bm{f}_S(\fnt{\tilde{u}}_{f,i}, \fnt{\tilde{u}}_{f,j})}\\
= &\frac{1}{2}\sum_{i\in I_J} A_i \bm{1}^Tn_{J,i}\LRp{ \sum_{j=1}^{|I_J|}c_{ij}\LRp{\psi(\fnt{\tilde{u}}_{f,j}) - \psi(\fnt{\tilde{u}}_{f,i})}}\\
= &\frac{1}{2} \bm{1}^Tn_{J,i}\LRp{ \sum_{i\in I_J} \sum_{j=1}^{|I_J|} A_ic_{ij}\LRp{\psi(\fnt{\tilde{u}}_{f,j}) - \psi(\fnt{\tilde{u}}_{f,i})}} \\
= &\frac{1}{2} \bm{1}^Tn_{J,i}\LRp{ \sum_{j=1}^{|I_J|} A_j \psi(\fnt{\tilde{u}}_{f,j}) - \sum_{i\in I_J} A_i\psi(\fnt{\tilde{u}}_{f,i})},
\end{align}
where we use the symmetry and conservation properties of the entropy conservative fluxes from Definition \ref{def:consevative_flux}. Distributing these contributions among all channels, the interface contributions on the $i$th channel are
\begin{align}
\frac{1}{2} \bm{1}^Tn_{J,i}\LRp{ \frac{1}{|I_J|}\sum_{j=1}^{|I_J|} A_j \psi(\fnt{\tilde{u}}_{f,j}) - A_i\psi(\fnt{\tilde{u}}_{f,i})}.
\label{eq:1d1d_ec_pf_flux_each_ch}
\end{align}
The flux contributions $\frac{1}{|I_J|}\sum_{j=1}^{|I_J|} A_j \psi(\fnt{\tilde{u}}_{f,j})$ in (\ref{eq:1d1d_ec_pf_flux_each_ch}) cancel at the junction interface because the neighboring elements have opposite normals. Combining the results at and away from the junction and using (\ref{eq:1D_1D_flux_coefficients}) with $\psi(\tilde{\fnt{u}}) - \fnt{v}^T \bm{f}^* = 0$ on non-junction boundaries, we reach the conclusion that
\begin{align}
\sum_{i\in I_J} A_i\sum_{k=1}^{K_i}\bm{1}^TJ^{k}_{i}\fnt{W}\frac{{\rm{d}}S(\fnt{u}_{q,i})}{{\rm{d}}t} = 0.
\end{align}
\end{proof}

We note that, for our numerical experiments, $c_{ij}$ are fixed values, but they can potentially vary with time or depend on solution values. However, determining expressions for junction coefficients $c_{ij}$ is beyond the scope of this paper. Fixed values of $c_{ij}$ still provide reasonable approximations in certain cases. We illustrate the application of this framework with two junction examples.
\subsubsection{Straight flow}
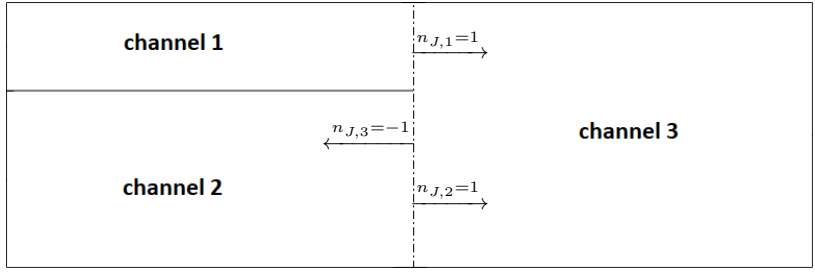
\begin{figure}[H]
\begin{center}
\begin{tikzpicture}
\node[anchor=south west,inner sep=0] at (0,0) {\includegraphics[width=0.65\textwidth]{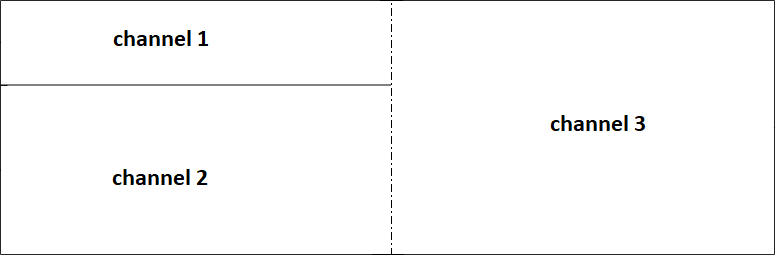}};
\node[font=\scriptsize] at (5.90,3) {$\xrightarrow{n_{J,1}=1}$};
\node[font=\scriptsize] at (5.90,1) {$\xrightarrow{n_{J,2}=1}$};
\node[font=\scriptsize] at (4.83,1.8) {$\xleftarrow{n_{J,3}=-1}$};
\end{tikzpicture}
\caption{1D-1D domain coupling.}
\label{fig:1_D_1D_mesh_split}
\end{center}
\end{figure}
First consider the following setting as in Figure \ref{fig:1_D_1D_mesh_split}. Let $A_1$,  $A_2$, and  $A_3$ be the width of the corresponding channels. We assume that $A_1 + A_2 = A_3$ and that the direction of the domains is aligned with the $x$ axis from left to right. Therefore, the outward surface normals are $n_1 = 1$, $n_2 = 1$, and $n_3 = -1$ for channels 1, 2, and 3 respectively as marked in Figure \ref{fig:1_D_1D_mesh_split}. We apply periodic boundary conditions in the $x$ direction. Since we are using an entropy stable DG scheme on each 1D domain, the only component left to specify is the treatment of the junction interface. Let $\bm{u}_1$, $\bm{u}_2$, and $\bm{u}_3$ be the conservative variables for channels 1, 2 and 3. We calculate the junction interface fluxes $\bm{f}_{J,i}$ in the following ways: for channels 1 and 2, we have
\begin{align}
\bm{f}_{J,1}(\bm{u}_1, \bm{u}_2, \bm{u}_3)= \bm{f}_S(\bm{u}_1, \bm{u}_3), \hspace{1cm} \bm{f}_{J,2}(\bm{u}_1, \bm{u}_2, \bm{u}_3)= \bm{f}_S(\bm{u}_2, \bm{u}_3).
\label{eq:1d_1d_sp}
\end{align}
When two channels join together, we consider the width of each channel when calculating flux contributions. Because channels 1 and 2 merge into channel 3, we average the flux contributions from channels 1 and 2 weighted by their widths:
\begin{align}
\bm{f}_{J,3}(\bm{u}_1, \bm{u}_2, \bm{u}_3)= \frac{\bm{f}_S(\bm{u}_1, \bm{u}_3)A_1 + \bm{f}_S(\bm{u}_2, \bm{u}_3)A_2}{A_3}.
\end{align}
This formulation falls under the general flux sharing framework in (\ref{eq:1D_1D_flux_coefficients}), with the following non-zero coefficients:
\begin{align}
c_{13} = c_{23} =1, \hspace{0.5cm} c_{31} = \frac{A_1}{A_3}, \hspace{0.5cm} c_{32} = \frac{A_2}{A_3}.
\end{align}
By Theorem \ref{thm:entropy_conservative_1D_1D}, this choice of coefficients is entropy conservative.
\subsubsection{Partial wall boundary condition}
\label{subsec:partial_wall}
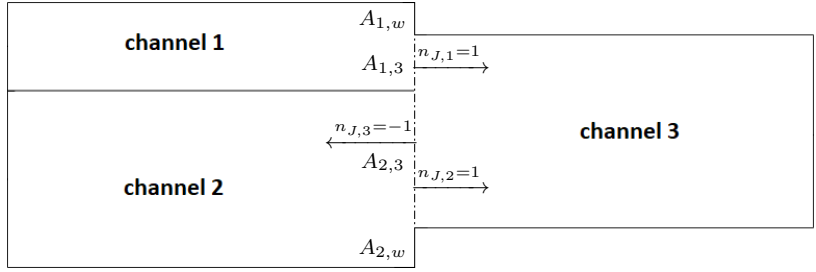
\begin{figure}[H]
\begin{center}
\begin{tikzpicture}
\node[anchor=south west,inner sep=0] at (0,0) {\includegraphics[width=0.65\textwidth]{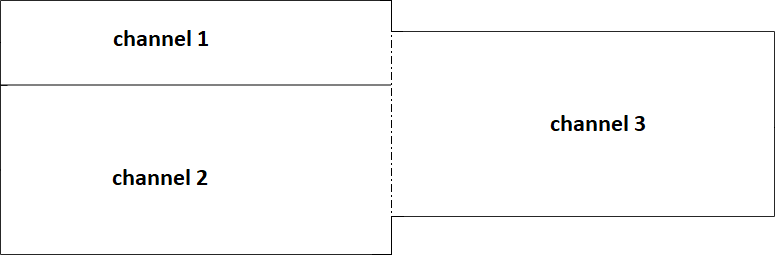}};
\node[font=\scriptsize] at (5,3.3) {$A_{1,w}$};
\node[font=\scriptsize] at (5,2.7) {$A_{1,3}$};
\node[font=\scriptsize] at (5,1.4) {$A_{2,3}$};
\node[font=\scriptsize] at (5,0.25) {$A_{2,w}$};
\node[font=\scriptsize] at (5.90,2.8) {$\xrightarrow{n_{J,1}=1}$};
\node[font=\scriptsize] at (5.90,1.2) {$\xrightarrow{n_{J,2}=1}$};
\node[font=\scriptsize] at (4.85,1.8) {$\xleftarrow{n_{J,3}=-1}$};
\end{tikzpicture}
\caption{Notation for 1D-1D junction coupling with ``partial'' wall boundary conditions.}
\label{fig:1_D_1D_mesh_ps_pwbc}
\end{center}
\end{figure}

In this second example, we introduce “partial” wall boundary conditions, where we blend the flux from a reflecting wall into the interface condition between channels 1, 2, and 3. We note that when the channel widths do not satisfy $A1 + A2 = A3$, the inclusion of a partial wall boundary condition is necessary to ensure entropy conservation. Without loss of generality, assume $A_1 + A_2 > A_3$ as shown in Figure \ref{fig:1_D_1D_mesh_ps_pwbc}. Let $A_{i,j}$ denote the width shared by channel $i$ and $j$ and let $A_{i,w}$ denote the width of the wall on channel $i$ at the junction. We then have:
\begin{align}
A_1 = A_{1,3}+A_{1,w}, \hspace{0.5cm} A_2 = A_{2,3}+A_{2,w}, \hspace{0.5cm} A_3 = A_{1,3}+A_{2,3}.
\end{align}
Again, recall that the direction of each domain is aligned with the $x$ axis from left to right such that the outward surface normals are $n_1 = 1$, $n_2 = 1$, and $n_3 = -1$ for channels 1, 2, and 3 respectively. Let $\bm{u}_i^{w}$denote the variables used in the entropy conservative flux to impose wall boundary conditions on the $i$th channel. For the shallow water equations, $\bm{u}_i^{w}$ is a ``mirror state'', such that
\begin{align}
\bm{u}_i^{w} = 
\begin{bmatrix}
h_i\\
-(hu)_i
\end{bmatrix},
\end{align}
where $h_i$ and $(hu)_i$ denote the water height and momentum on the $i$th channel \cite{wintermeyer2017entropy, chen2017entropy}. We calculate the fluxes in the following ways:
\begin{align}
\bm{f}_{J,1}(\bm{u}_1, \bm{u}_1^{w}, \bm{u}_2,\bm{u}_2^{w}, \bm{u}_3)&= \frac{\bm{f}_S(\bm{u}_1, \bm{u}_3)A_{1,3} + \bm{f}_S(\bm{u}_1, \bm{u}_1^{w})A_{1,w}}{A_1}, \\
\bm{f}_{J,2}(\bm{u}_1, \bm{u}_1^{w}, \bm{u}_2,\bm{u}_2^{w}, \bm{u}_3)&= \frac{\bm{f}_S(\bm{u}_2, \bm{u}_3)A_{2,3} + \bm{f}_S(\bm{u}_2, \bm{u}_2^{w})A_{2,w}}{A_2}, \\
\bm{f}_{J,3}(\bm{u}_1, \bm{u}_1^{w}, \bm{u}_2,\bm{u}_2^{w}, \bm{u}_3)&= \frac{\bm{f}_S(\bm{u}_1, \bm{u}_3)A_{1,3} + \bm{f}_S(\bm{u}_2, \bm{u}_3)A_{2,3}}{A_3}.
\label{eq:pwbc_flux}
\end{align}
For the converged channel 3 in (\ref{eq:pwbc_flux}), we average the flux from channels 1 and 2 weighted by the width of the part of the channel that are shared with channel 3. However, for channels 1 and 2, we average the fluxes to capture interactions with channel 3 as well as from the imposition of wall boundary conditions weighted by their respective widths.

The ``partial'' wall boundary condition provides a treatment for channel junctions with non-matching widths, which preserves entropy conservation under entropy conservative wall boundary conditions \cite{chen2017entropy}. The partial wall boundary conditions correspond to the general junction treatment (\ref{eq:1D_1D_flux_coefficients}) with the following coefficients:
\begin{align}
\nonumber
c_{1,1} &= \frac{A_{1,w}}{A_1}, &c_{1,2} &= 0, &c_{1,3} &= \frac{A_{1,3}}{A_1}\\
\nonumber
c_{2,1} &= 0, &c_{2,2} &=\frac{A_{2,w}}{A_2}, &c_{2,3} &=\frac{A_{2,3}}{A_2}  \\
\nonumber
c_{3,1} &= \frac{A_{1,3}}{A_3}, &c_{3,2} &= \frac{A_{2,3}}{A_3}, &c_{3,3} &=0.
\end{align}
Here, note that the coefficients $c_{11}$ and $c_{22}$ scale $\bm{f}_S(\bm{u}_i, \bm{u}_i^{w})$ rather than $\bm{f}_S(\bm{u}_i, \bm{u}_i)$.
\subsection{Entropy dissipation at interfaces}
The final step of building the entropy stable DG method is adding entropy dissipation terms at the coupling interfaces. To accomplish this, we apply Lax-Friedrichs penalization \cite{chan2018discretely}. Local Lax-Friedrichs penalization augments the flux function at element interfaces with an additional term:
\begin{align}
\bm{f}_{S}(\bm{u}, \bm{u}^+) \xrightarrow{} \bm{f}_{S}(\bm{u}, \bm{u}^+) - \frac{\lambda}{2} \llbracket\bm{u}\rrbracket \quad \textrm{in 1D},\\
\bm{f}_{i,S}(\bm{u}, \bm{u}^+) \xrightarrow{} \bm{f}_{i,S}(\bm{u}, \bm{u}^+) - \frac{\lambda}{2} \llbracket\bm{u}\rrbracket \quad \textrm{in 2D},
\end{align}
where $\lambda$ is an estimate of the maximum wave speed and $\llbracket\cdot\rrbracket$ denotes the jump across the interface, $\llbracket\bm{u}\rrbracket = \bm{u}^+ -\bm{u}$. While defining $\llbracket \bm{u}\rrbracket$ is straightforward on meshes in a single dimension, the procedure is more involved for 1D-2D meshes. Let $\bm{f}_{J,2D}^p$ denote the penalized 2D flux. We first compute the flux as specified in (\ref{eq:1D-2D_flux2D}), then adding the Lax-Friedrichs penalization yields
\begin{align}
\bm{f}_{J,2D}^p(\bm{u}_{1D}, \bm{u}_{2D})  = \sum_{i=1}^d n_i\bm{f}_{i,S}(\bm{u}_{2D}, \bm{u}_{2D}^+) - \frac{\lambda}{2} \llbracket\bm{u}_{2D}\rrbracket,\qquad \bm{u}_{2D}^+ = \bm{R}\bm{u}_{1D},
\end{align}
where $n_i$ is the $i$th component of the unit outward normal on the 2D side of the 1D-2D interface. 
Then, we can define the 1D flux at the interface in a manner similar to equation (\ref{eq:1D-2D_flux}):
\begin{align}
\bm{f}_{J,1D,p}(\bm{u}_{1D}, \bm{u}_{2D}) = \frac{\fnt{w}_{J,f}^T\LRp{\bm{R}^T\bm{f}_{J,2D}^p(\bm{u}_{1D}, \bm{u}_{2D})}}{\fnt{w}_{J,f}^T\bm{1}}.
\label{eq:1D-2D_flux_LF}
\end{align}
We pass the penalized fluxes calculated at the 1D-2D interface to the 1D domain. We then scale the fluxes in 1D by the surface normal to ensure they have the correct sign. For mesh elements not on the 1D-2D interface, we use standard 1D or 2D Lax-Friedrichs penalization in their respective domains \cite{chen2017entropy}.

For 1D-1D junction treatments, to implement Lax-Friedrichs penalization, we simply make the following changes to fluxes between solutions on each channel:
\begin{align}
\bm{f}_{S}(\bm{u}_{f,i}, \bm{u}_{f,j}) \xrightarrow{} \bm{f}_{S}(\bm{u}_{f,i}, \bm{u}_{f,j}) - \frac{\lambda}{2} \llbracket\bm{u}\rrbracket.
\end{align}
\end{subequations}

%% file: section5.tex
\section{Numerical results}
\label{sec:numerical_results}
\begin{subequations}
\numberwithin{equation}{section}
In this section, we present numerical experiments to demonstrate the accuracy and stability of the entropy stable DG scheme with 1D-2D and 1D-1D junction couplings. The first experiment is a simple parallel ``split and converge'' channel network as shown in Figure \ref{fig:2D_mesh_split}. In the second experiment, the setup consists of a channel ``split and converge'' in the shape of diamond, as shown in Figure \ref{fig:diamond_mesh}. We test both 1D-2D and 1D-1D coupling on these two setups. The third experiment studies the T-shaped junction as shown in Figure \ref{fig:2D_mesh_T}. The fourth experiment focuses on channels turning at different angles as shown in Figure \ref{fig:2D_mesh_turns}. In the last experiment, we simulate a dam break on the domain shown in Figure \ref{fig:Dam_break_mesh}. Here, a large reservoir is connected to a long channel with a $45^\circ$ turn in the middle. 

We perform computations with the shallow water equations on all setups. We also present experiments for the compressible Euler equations on select setups in Appendix A. For the shallow water equations, we have the following well-balanced and entropy conservative fluxes in 2D \cite{tadmor2003entropy, carpenter2014entropy, gassner2016well, wintermeyer2018entropy}
\begin{align}
\bm{f}^x_{S}\LRp{\bm{u}_L,\bm{u}_R} &=
\begin{bmatrix}
\avg{hu}\\
\avg{hu}\avg{u} + g\avg{h}^2 - \frac{1}{2}g\avg{h^2}\\
\avg{hu}\avg{v}
\end{bmatrix}
\label{eq:SWE_flux2dx}
, \\
\bm{f}^y_{S}\LRp{\bm{u}_L,\bm{u}_R} &=
\begin{bmatrix}
\avg{hv}\\
\avg{hv}\avg{u}\\
\avg{hv}\avg{v} + g\avg{h}^2 - \frac{1}{2}g\avg{h^2}\\ 
\end{bmatrix}.
\label{eq:SWE_flux2dy}
\end{align}
The shallow water equations fluxes in 1D are
\begin{align}
\bm{f}_{S1D}\LRp{\bm{u}_L,\bm{u}_R} &=
\begin{bmatrix}
\avg{hu}\\
\avg{hu}\avg{u} + \frac{1}{2}gh_Lh_R\\
\end{bmatrix}.
\label{eq:SWE_flux1D}
\end{align}

For all experiments, we use 4th order 5-stage low storage Runge-Kutta method of \cite{carpenter1994fourth}. Following the derivation of stable timestep restrictions in \cite{chan2016gpu}, we define the timestep $\Delta t$ as the following:
\begin{align}
&C_{N2D} = \frac{(N_{2D}+1)(N_{2D}+2)}{2},\hspace{1cm} C_{N1D} = \frac{(N_{1D}+1)^2}{2},\\
&\Delta t = min\LRp{CFL \times \frac{h_{1D}}{C_{N1D}}, CFL \times \frac{h_{2D}}{C_{N2D}}},
\end{align}
where $C_N$ is the degree dependent constant in the inverse trace inequality \cite{warburton2003constants}, and CFL is a user-defined constant. We use $CFL = 0.25$ and $g = 1$ for all experiments. 
When comparing 1D and 2D solutions, we average the 2D solution along the width of the channel.

\subsection{Shallow water equation experiments}
\subsubsection{Parallel split and converge (1D-2D and 1D-1D junction treatments)}
\label{subsec:ps_swe}
We first consider the parallel ``split and converge'' geometry shown in Figure \ref{fig:2D_mesh_split}. We have a channel which split into two parallel channels. The fully 2D domain ranges from $-4$ to $4$ in $x$ direction and $-1$ to $1$ in $y$ direction. We apply periodic boundary conditions in the $x$ direction and wall boundary conditions in $y$ direction for the fullly 2D simulation. We consider both fully 2D, 1D-2D coupling and 1D-1D junction coupling for this case. A diagram of the fully 2D and 1D-2D couplings are shown in Figure \ref{fig:2D_mesh_split}. There are 32 uniform elements of size 0.25 in all 1D domains and 128 uniform triangles for the fully 2D mesh. We have parallel channels separated by a wall marked in a red dashed line on the mesh in Figure \ref{fig:2D_mesh_split}.
\begin{figure}
\begin{center}
\includegraphics[width=.45\textwidth]{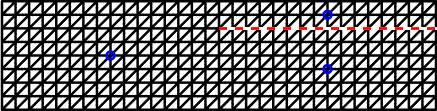}
\qquad
\includegraphics[width=.45\textwidth]{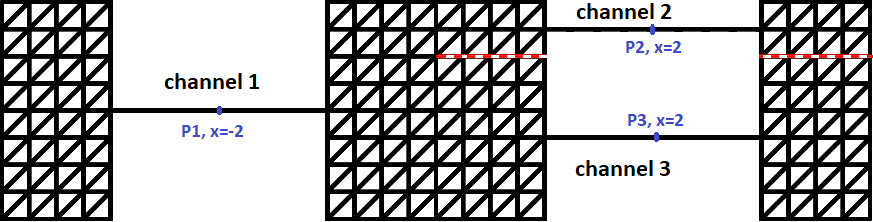}
\caption{Fully 2D mesh and 1D-2D coupling mesh of parallel split and converge with wall marked in red dash line.}
\label{fig:2D_mesh_split}
\end{center}
\end{figure}
We use the following discontinuous initial conditions:
\begin{align}
h_0 = \begin{cases}
      3 & \text{in channel 1}\\
      4 & \text{in channel 2 and 3}\\
    \end{cases}, 
\hspace{1cm} u_0=v_0=0.
\label{eq:SWE_ps_initial_condition_1}
\end{align}
We run this experiment until final time $T=2$. We have the semi-discrete system 
\begin{align}
M\td{\bm{u}}{t} = r(\bm{u}).
\label{eq:entropy_rhs}
\end{align}
To verify entropy conservation, we compute the entropy right hand side (RHS) $\bm{v}^Tr(u)$, without Lax-Friedrichs dissipation. Note that, according to Theorem \ref{thm:entropy_conservative} and \ref{thm:entropy_conservative_1D_1D}, $\bm{v}^Tr(u)$ should be zero in the absence of Lax-Friedrichs dissipation.
\begin{table}[H]
\begin{center}
\begin{tabular}{|c|c|c|c|}
\hline
                 &  $N_{2D} = 3$  &  $N_{2D} = 4$  &  $N_{2D} = 5$ \\
\hline
   $N_{1D} = 3$  &  2.6468e-13 & 1.0147e-13 & 4.4698e-13 \\
\hline
   $N_{1D} = 4$  &  2.0872e-13 & 1.0147e-13 & 4.6718e-13 \\
\hline
   $N_{1D} = 5$  &  7.7982e-13 & 9.8765e-13 & 7.8337e-13 \\
\hline
\end{tabular}
\end{center}
\caption{Maximum of absolute value of entropy RHS for SWE 1D-2D coupling.}
\label{Tab:swe_entopy_rhs}
\end{table}
We test this setting without Lax-Friedrichs dissipation, We also  utilize different polynomial degrees $N_{1D}$ and $N_{2D}$ in each domain. The maximum absolute values of the entropy RHS (\ref{eq:entropy_rhs}) throughout the run are listed in Table \ref{Tab:swe_entopy_rhs}. We also list the values of the entropy RHS for the 1D-1D coupling in Table \ref{Tab:1D_1D_entopy_rhs}. We observe that the entropy RHS values are always close to machine precision, which confirms that our coupling scheme is entropy conservative for the shallow water equations. 

To test accuracy, we compare the values of the conservative variables at the midpoints of each channel, which are marked in blue on the right side in Figure \ref{fig:2D_mesh_split}. We plot the water height at each midpoint for the fully 2D, 1D-2D, and 1D-1D junction treatment and compare the results in Figure \ref{fig:SWE_parallel_split_points_time}. We set the polynomial degree to $N=3$ in all domains and use the local Lax-Friedrichs penalization to add entropy dissipation. Note that oscillation appear in the numerical solution. This is due to the presence of shock discontinuities and absence of slope limiters and artificial viscosity in our numerical scheme. However, the solutions remain stable and do not blow up due to entropy stability.
\begin{figure}
\centering
\subfloat[P1]{\includegraphics[width=.32\textwidth]{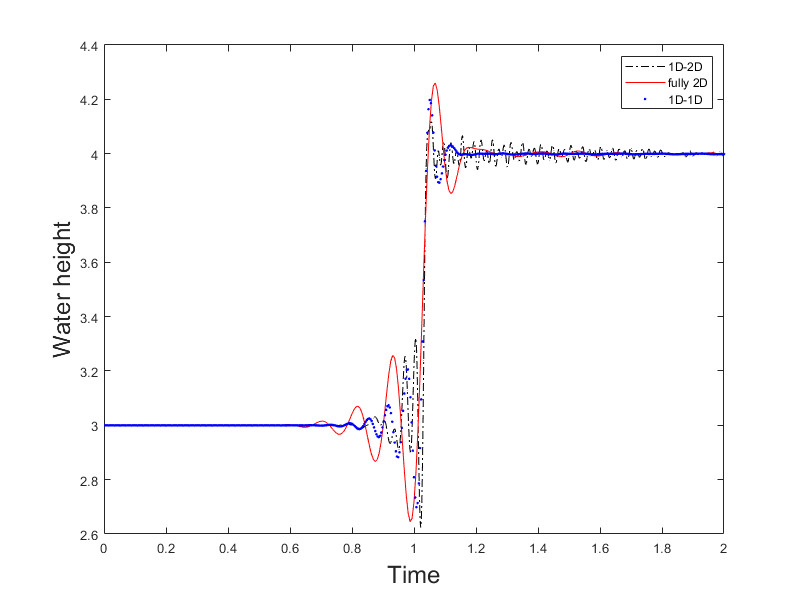}}
\hspace{.1em}
\subfloat[P2]{\includegraphics[width=.32\textwidth]{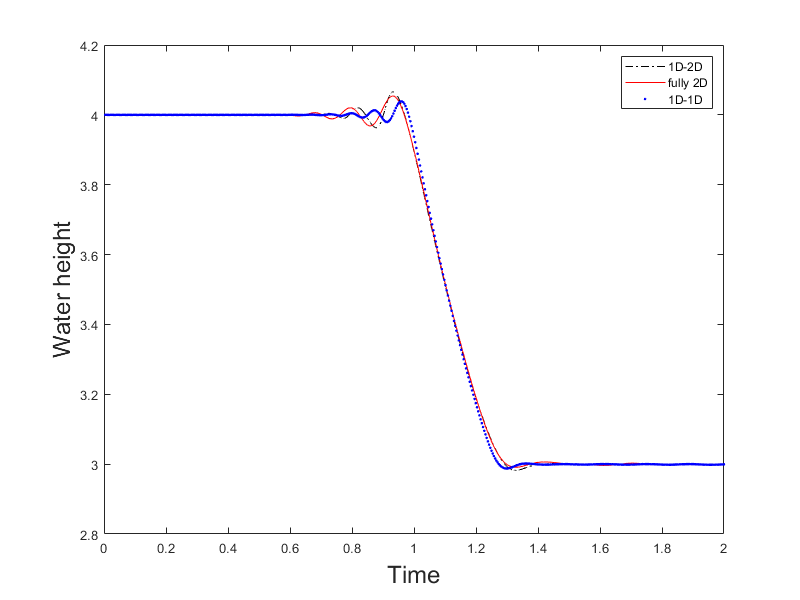}}
\hspace{.1em}
\subfloat[P3]{\includegraphics[width=.32\textwidth]{figs/swe_parallel_split_p2_3_2s.png}}
\caption{Water height in parallel split with initial conditions (\ref{eq:SWE_ps_initial_condition_1}) at points P1 (a), P2 (b), and P3 (c) for the fully 2D, 1D-2D, and 1D-1D junction models.}
\label{fig:SWE_parallel_split_points_time}
\end{figure}

We also test with the following initial conditions:
\begin{align}
h_0 = \begin{cases}
      4 & \text{in channel 1}\\
      5 & \text{in channel 2}\\
      6 & \text{in channel 3}
    \end{cases}, 
\hspace{1cm} u_0=v_0=0.
\label{eq:SWE_ps_initial_condition_2}
\end{align}
We plot the water height at points P1, P2, and P3 under initial conditions (\ref{eq:SWE_ps_initial_condition_2}) in Figure \ref{fig:SWE_parallel_split_points_time_456}.
\begin{figure}
\centering
\subfloat[P1]{\includegraphics[width=.32\textwidth]{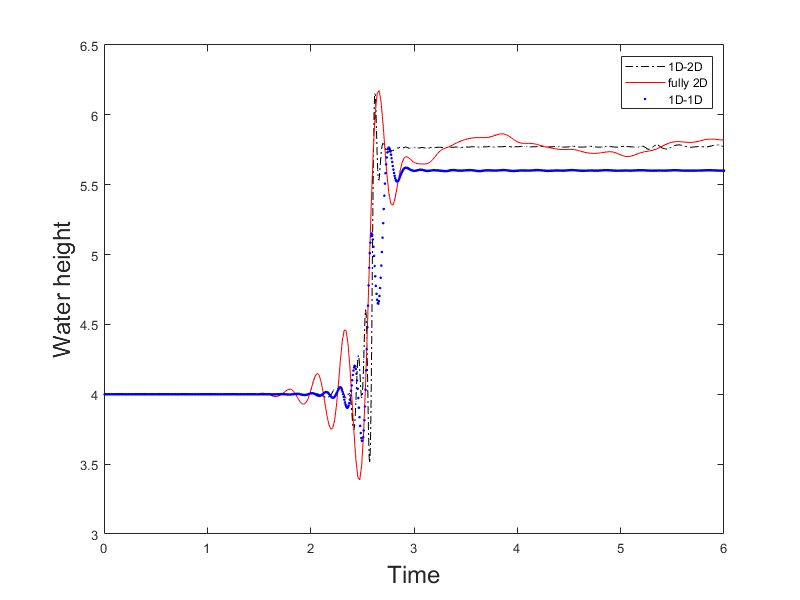}}
\hspace{.1em}
\subfloat[P2]{\includegraphics[width=.32\textwidth]{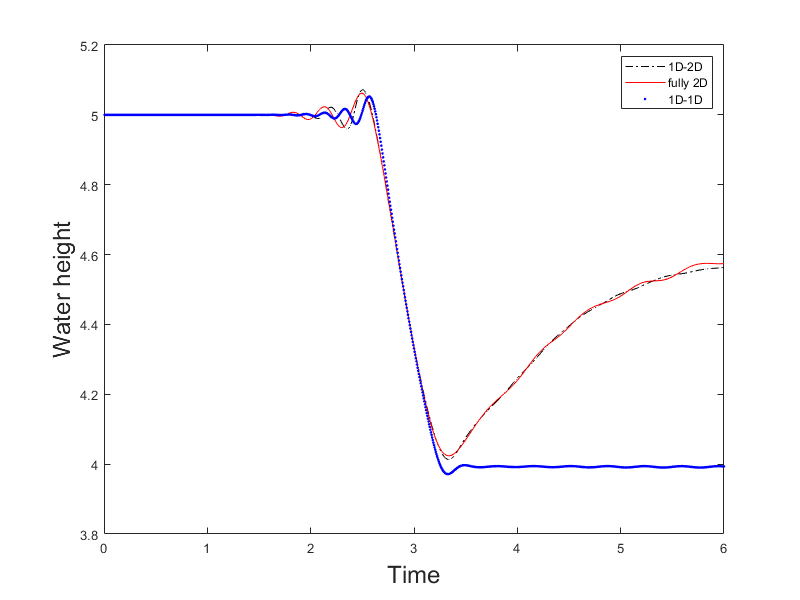}}
\hspace{.1em}
\subfloat[P3]{\includegraphics[width=.32\textwidth]{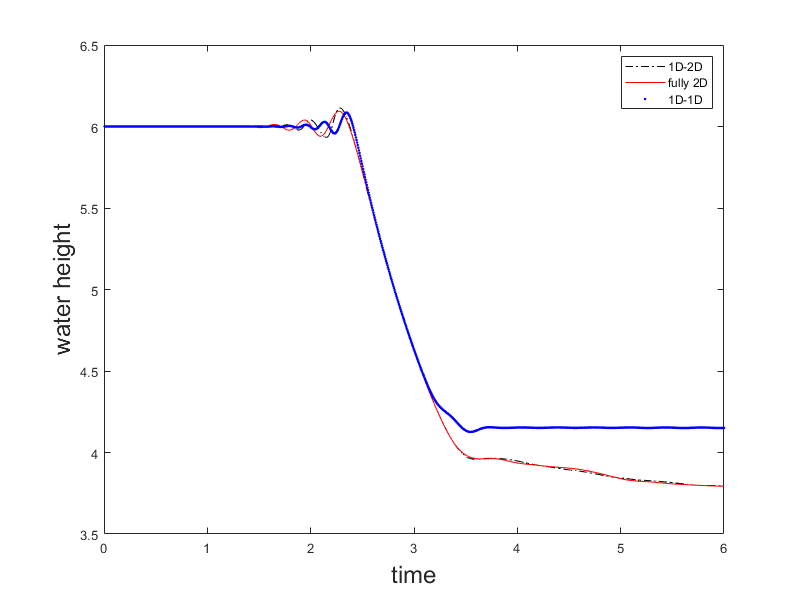}}
\caption{Water height in parallel split with initial conditions (\ref{eq:SWE_ps_initial_condition_2}) at points P1 (a), P2 (b), and P3 (c) for the fully 2D, 1D-2D, and 1D-1D junction models.}
\label{fig:SWE_parallel_split_points_time_456}
\end{figure}
We notice that the reduced 1D-2D model produces results close to the fully 2D model as was observed in \cite{bellamoli2018numerical, neupane2015discontinuous}. The 1D-1D model deviates from the 2D model more significantly under initial conditions (\ref{eq:SWE_ps_initial_condition_2}).  With initial conditions (\ref{eq:SWE_ps_initial_condition_2}),  there is significant vertical water motion at the junction which the 1D-1D model fails to capture.

We also consider the parallel split with non-matching widths, where the sum of the widths of channels 2 and 3 does not match the width of channel 1. The mesh is shown in Figure \ref{fig:SWE_psdw_mesh}, where channels 1, 2, and 3 have widths of $\sqrt{2}$, 1, and 1. Each channel has a length of 4 for the fully 2D simulation. We enforce wall boundary conditions in both the $x$ and $y$ direction. We use the mesh generator within PDEModel from MATLAB to construct our 2D mesh \cite{MatlabPDETB}, with $h_{min} = 0.25$. In the 1D-1D model, each 1D domain consist of 32 uniform elements of size 0.25.
\begin{figure}
\centering
\subfloat[]{\raisebox{1ex}{\includegraphics[width=.45\textwidth]{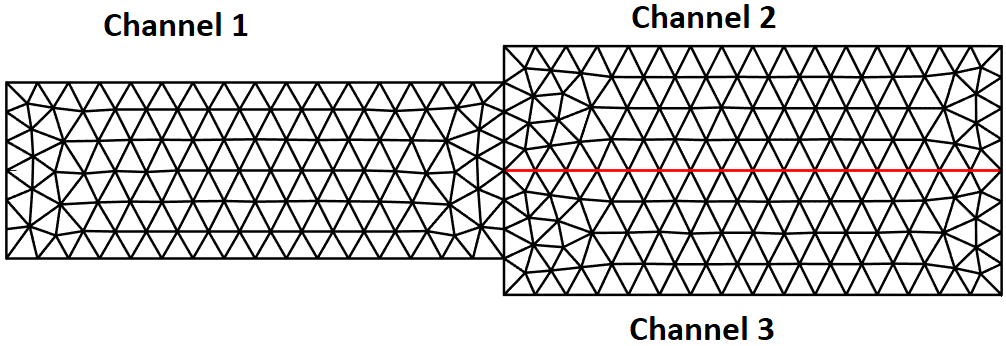}}}\qquad
\subfloat[]{\includegraphics[width=.45\textwidth]{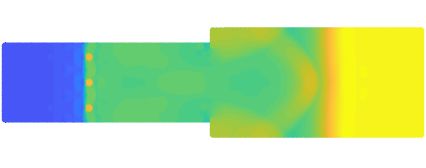}}
\caption{Parallel split with change in width at junction.}
\label{fig:SWE_psdw_mesh}
\end{figure}

For 1D-1D model, we implement partial wall boundary conditions at the junction. We use the initial conditions (\ref{eq:SWE_ps_initial_condition_1}) and plot the water height at the midpoints P1, P2, and P3 from the fully 2D model and from 1D-1D model in Figure \ref{fig:SWE_psdw_points}. We observe that the 1D model performs reasonably well under such conditions.
\begin{figure}[H]
\centering
\subfloat[P1]{\includegraphics[width=.32\textwidth]{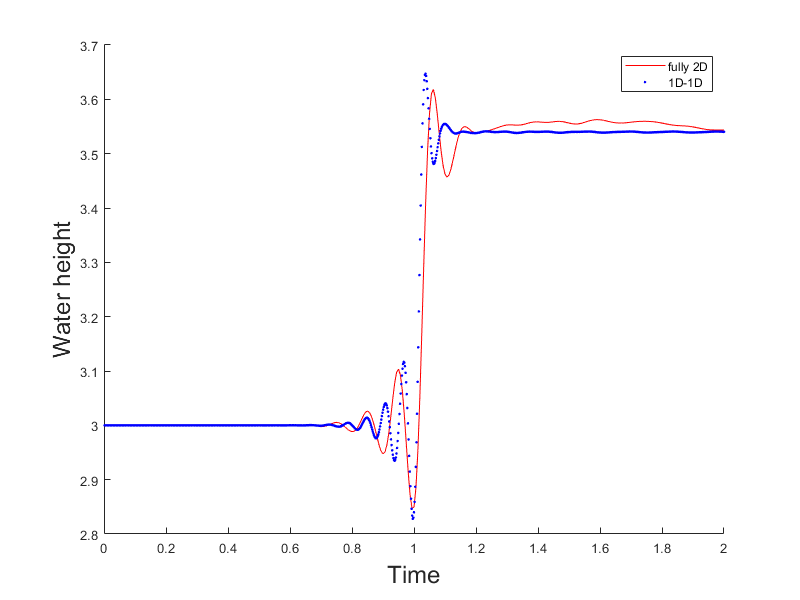}}
\hspace{.1em}
\subfloat[P2]{\includegraphics[width=.32\textwidth]{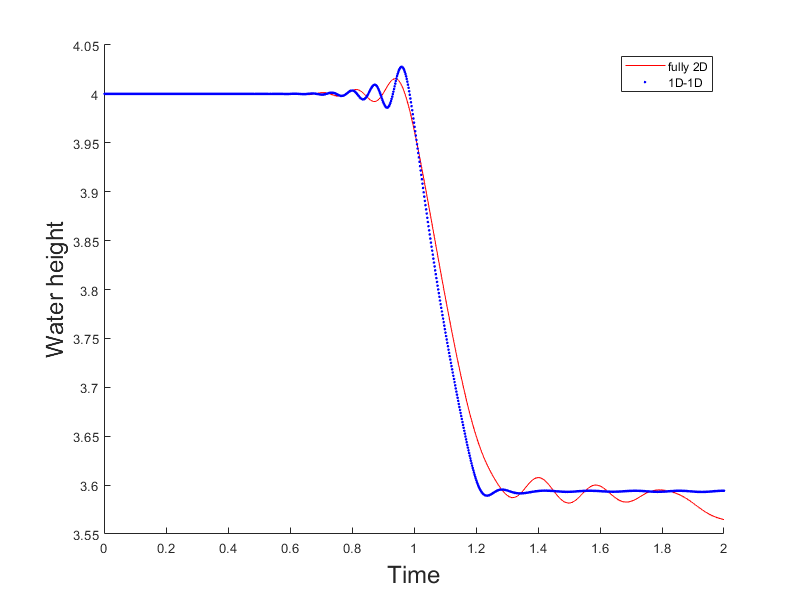}}
\hspace{.1em}
\subfloat[P3]{\includegraphics[width=.32\textwidth]{figs/swe_psdf_P2.png}}
\caption{Water height in parallel split with non-matching widths at junction under initial conditions (\ref{eq:SWE_ps_initial_condition_2}) at points P1 (a), P2 (b), and P3 (c) for the fully 2D, and 1D-1D junction models.}
\label{fig:SWE_psdw_points}
\end{figure}

\subsubsection{Diamond split and converge (1D-2D and 1D-1D junction treatments)}
\label{subsec:diamond_swe}
\begin{figure}[H]
\centering
\subfloat[]{\includegraphics[height=.15\textheight]{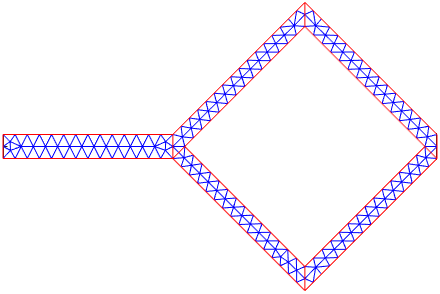}}
\subfloat[]{\includegraphics[height=.15\textheight]{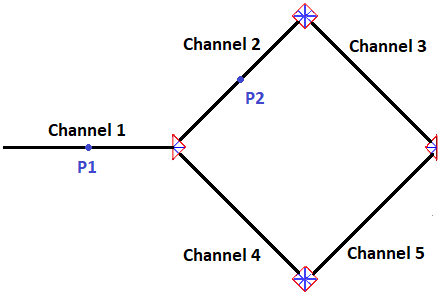}}
\subfloat[]{\raisebox{-4.5ex}{\includegraphics[height=.22\textheight]{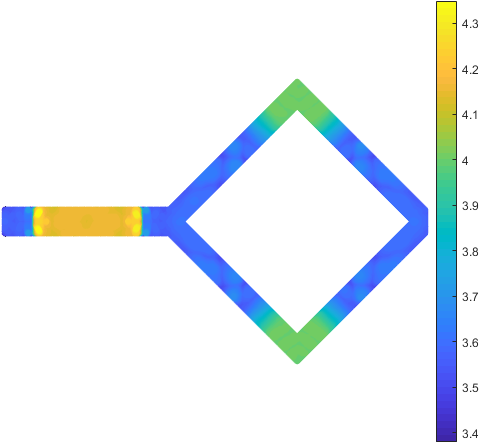}}}
\caption{Diamond split and converge. Fully 2D mesh (a), 1D-2D mesh (b), and solution snapshot (c).}
\label{fig:diamond_mesh}
\end{figure}
The second numerical experiment consists of a diamond split and converge setup, as shown in Figure \ref{fig:diamond_mesh}. We enforce wall boundary conditions everywhere except at the left side of the channel 1, which connects back to the right side of the domain to enforce periodicity. The horizontal channel has a width of $\sqrt{2}$ and each channel in the diamond has a width of $1$. The length of each 1D channel is 10, not including the junction. We visualize the solution at $P1$ and $P2$, the midpoints of channels 1 and 2. We use the mesh generator within PDEModel from MATLAB to construct our 2D mesh \cite{MatlabPDETB}, with $h_{min} = \sqrt{2}$. 

The 1D-2D junction treatment is shown in Figure \ref{fig:diamond_mesh}, where we represent junction using triangular meshes. Each 1D domain consists of 16 uniform elements of size 0.625. For this diamond split, we also implement a 1D-1D coupling scheme with partial wall boundary conditions at the junctions. We note that the 1D-1D junction treatments do not account for channel angles.

The initial conditions for shallow water equations with this setup are 
\begin{align}
h_0 = \begin{cases}
      3 & \text{in channel 1}\\
      4 & \text{otherwise}\\
    \end{cases}, 
\hspace{1cm} u_0=v_0=0.
\label{eq:SWE_diamond_initial_condition}
\end{align}
\begin{figure}[H]
\centering
\subfloat[P1]{\includegraphics[width=.4\textwidth]{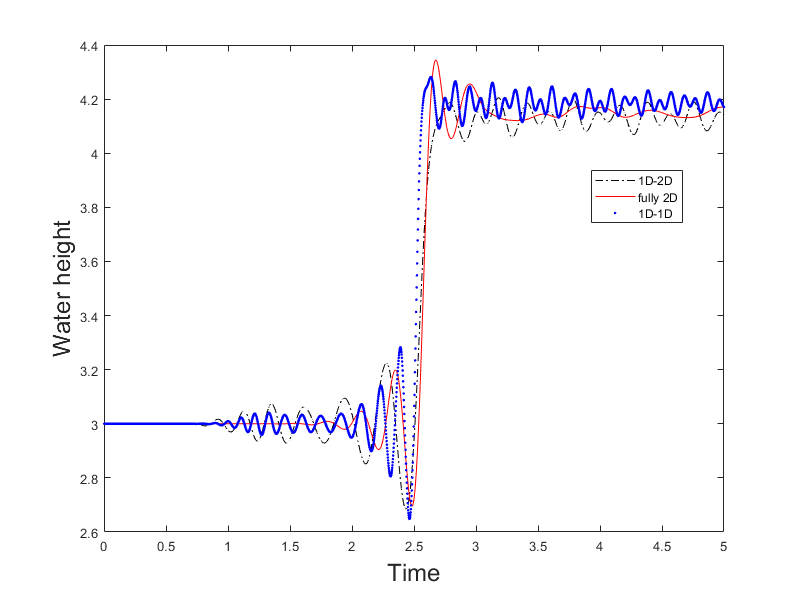}}
\hspace{1em}
\subfloat[P2]{\includegraphics[width=.4\textwidth]{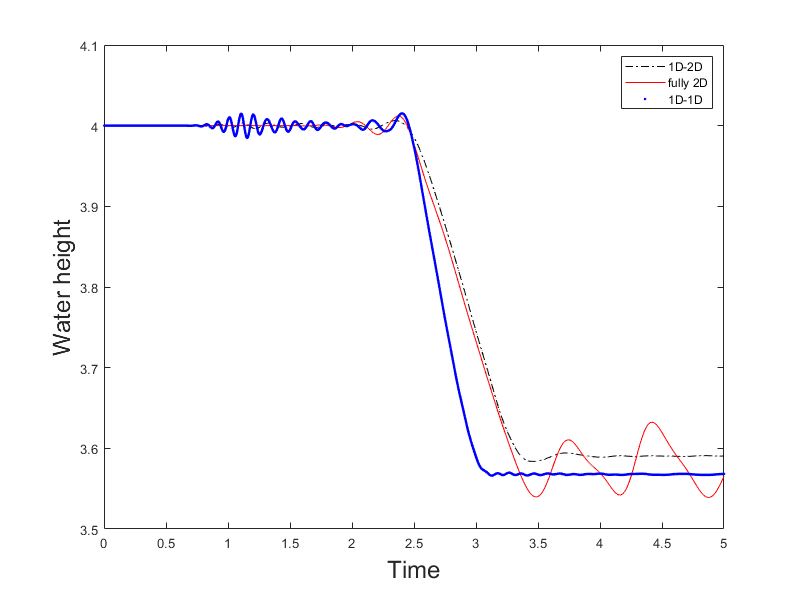}}
\caption{Water height in diamond split with initial conditions (\ref{eq:SWE_diamond_initial_condition}) at points P1 (a) and P2 (b) for the fully 2D, 1D-2D, and 1D-1D junction models.}
\label{fig:diamond_swe_p}
\end{figure}
The computed entropy RHS (\ref{eq:entropy_rhs}) during the run are on the order of $10^{-15}$ to $10^{-13}$ when not applying Lax-Friedrichs dissipation, which verifies entropy conservation. To demonstrate the accuracy of the 1D-2D junction treatment, we use the same initial conditions with Lax-Friedrichs dissipation. We plot the water height at points $P1$ and $P2$ in Figure \ref{fig:diamond_swe_p}. We observe that the solutions from the 1D-2D model match the arrival time of the solutions from the fully 2D model, but display discrepancies in amplitude at $P2$. The results from 1D-1D model produce an earlier arrival time for the shock at $P2$. 

\subsubsection{T-junction (1D-1D junction treatment)}
For the T-junction experiment, we present the setup in Figure \ref{fig:2D_mesh_T}. Since previous experiments and other authors \cite{bellamoli2018numerical, neupane2015discontinuous} confirmed the accuracy of the 1D-2D junction treatments, we focus on the comparison between the 1D-1D junction model and the fully 2D junction model for this case. Each channel has a width of 1 and length of 10. We define $P1$, $P2$, and $P3$ as the midpoints of each channel and record the water heights at these points. PDEModel is used to construct the fully 2D mesh, with $h_{min} = 0.25$. For 1D-1D model, we represent each channel in 32 uniform elements of size 0.3125. We use the flux sharing framework with the following coefficients:
\begin{align}
c_{ij} = \begin{cases}
      0  & \text{$i = j$}\\
      \frac{1}{2}  & \text{$i\neq j$}\\
    \end{cases}.
\end{align} 
We use polynomial degree $N=3$ in all experiments and run up to time $T=6$. We turn on Lax-Friedrichs dissipation to compare solutions at midpoints of each 1D channel with solutions from the fully 2D model. We test on the following three sets of initial conditions for the shallow water equations:
\begin{figure}
\begin{center}
\includegraphics[width=.16\textheight]{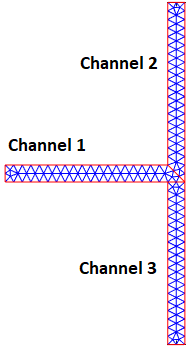} \hspace{2cm}
\includegraphics[width=.16\textheight]{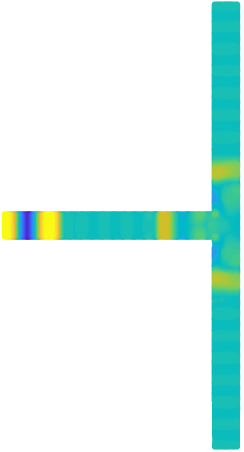}
\caption{Fully 2D mesh and snapshot of solution on T-junction.}
\label{fig:2D_mesh_T}
\end{center}
\end{figure}

\begin{align}
h_0 = \begin{cases}
      6  & \text{in $x\leq4$ channel 1}\\
      4  & \text{otherwise}\\
    \end{cases}, 
\hspace{1cm} u_0=v_0=0.
\label{eq:SWE_T_initial_condition_1}
\end{align}
\begin{align}
h_0 = \begin{cases}
      6  & \text{in $x\leq4$ channel 1}\\
      6  & \text{in $y\geq4$ channel 2}\\
      4  & \text{otherwise}\\
    \end{cases}, 
\hspace{1cm} u_0=v_0=0.
\label{eq:SWE_T_initial_condition_2}
\end{align}
\begin{align}
h_0 = \begin{cases}
      6  & \text{in $x\leq4$ channel 1}\\
      5  & \text{in $y\geq4$ channel 2}\\
      5.5  & \text{in $y\leq-4$ channel 3}\\
      4  & \text{otherwise}\\
    \end{cases}, 
\hspace{1cm} u_0=v_0=0.
\label{eq:SWE_T_initial_condition_3}
\end{align}
These initial conditions correspond to shocks propagating from one, two, and all three channels. We plot the water heights at $P1$, $P2$, and $P3$ in Figure \ref{fig:SWE_T_points_1}, \ref{fig:SWE_T_points_2}, and \ref{fig:SWE_T_points_3} respectively.
\begin{figure}
\centering
\subfloat[P1]{\includegraphics[width=.32\textwidth]{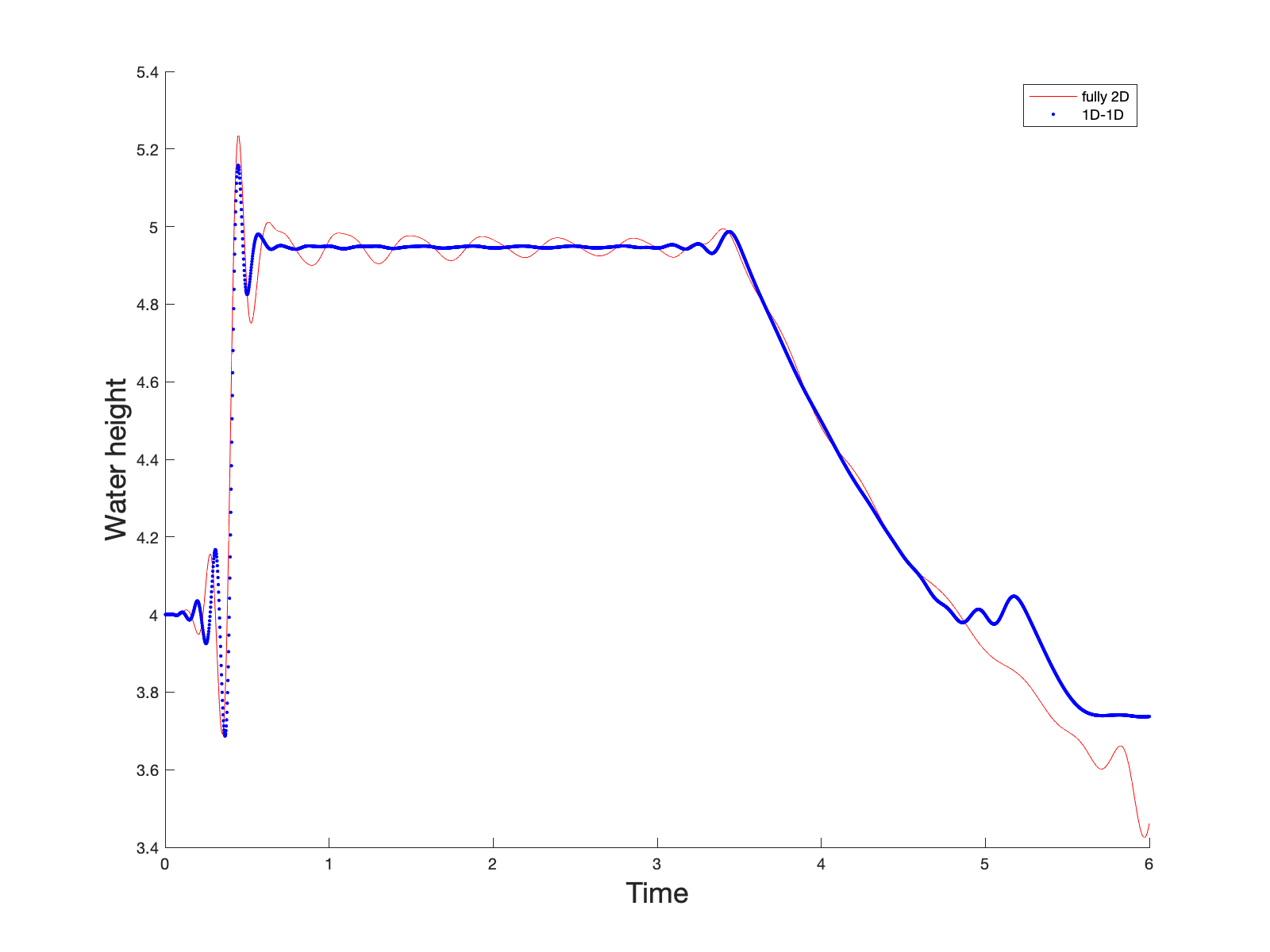}}
\hspace{.1em}
\subfloat[P2]{\includegraphics[width=.32\textwidth]{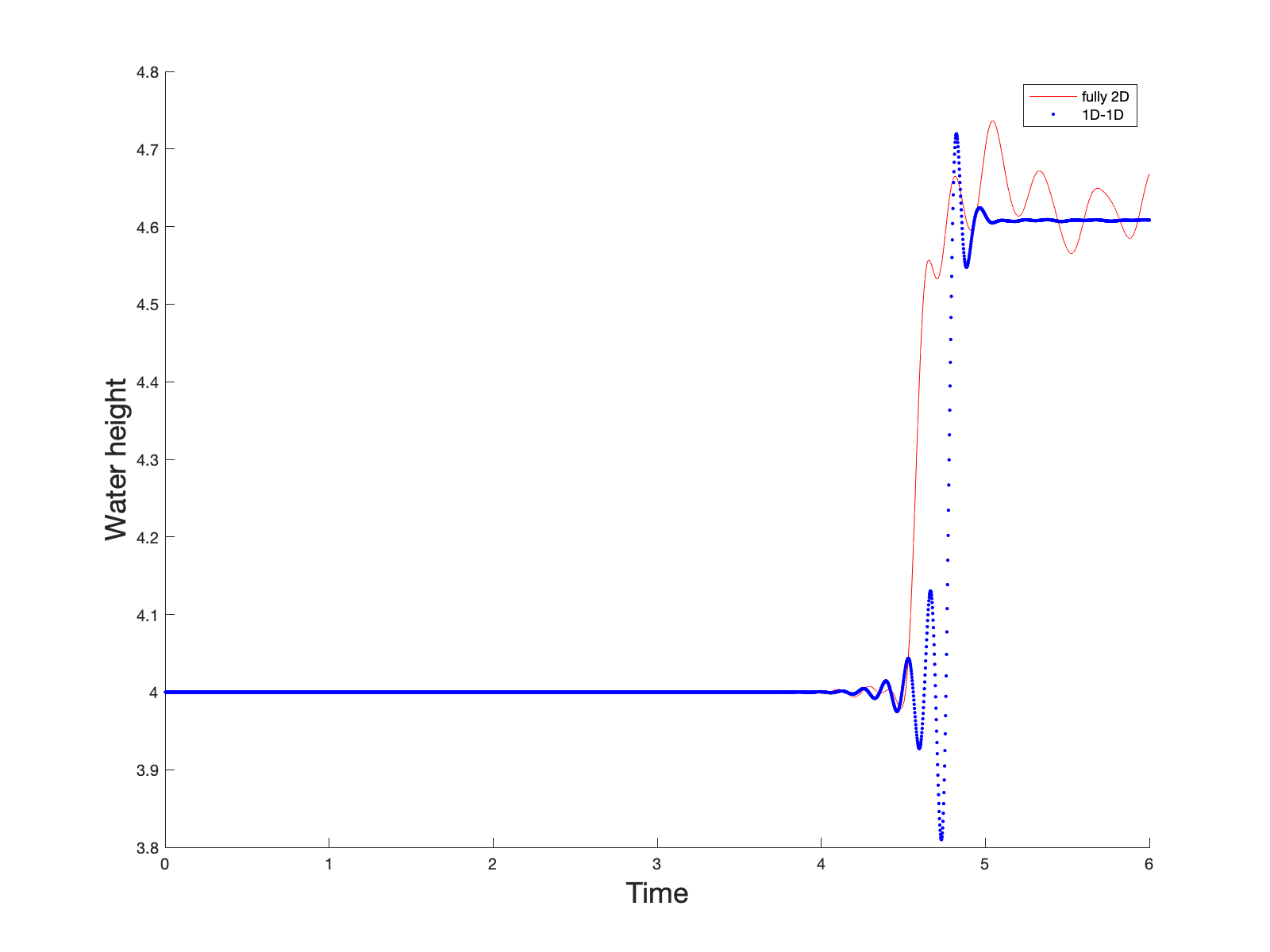}}
\hspace{.1em}
\subfloat[P3]{\includegraphics[width=.32\textwidth]{figs/swe_T_P2.png}}
\caption{Water height in T-junction with initial conditions (\ref{eq:SWE_T_initial_condition_1}) at points P1 (a), P2 (b), and P3 (c) for the fully 2D and 1D-1D junction models.}
\label{fig:SWE_T_points_1}
\end{figure}
\begin{figure}
\centering
\subfloat[P1]{\includegraphics[width=.32\textwidth]{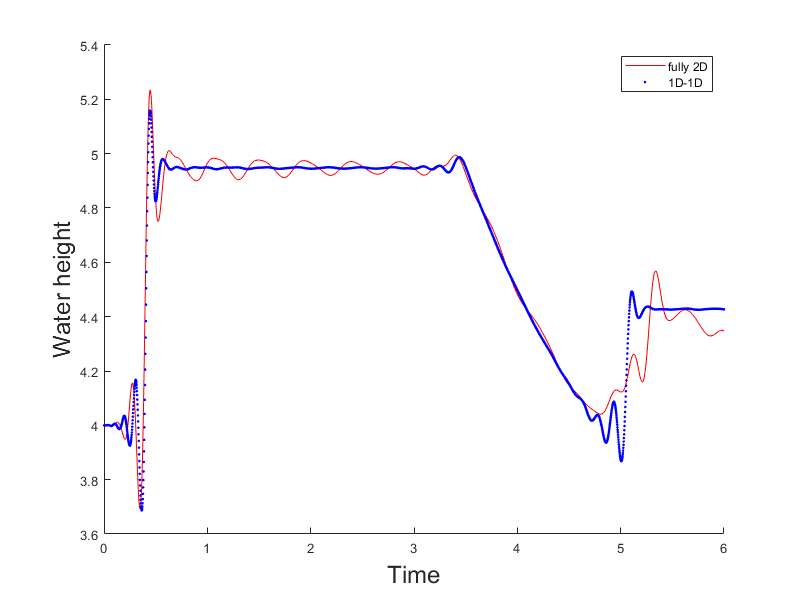}}
\hspace{.1em}
\subfloat[P2]{\includegraphics[width=.32\textwidth]{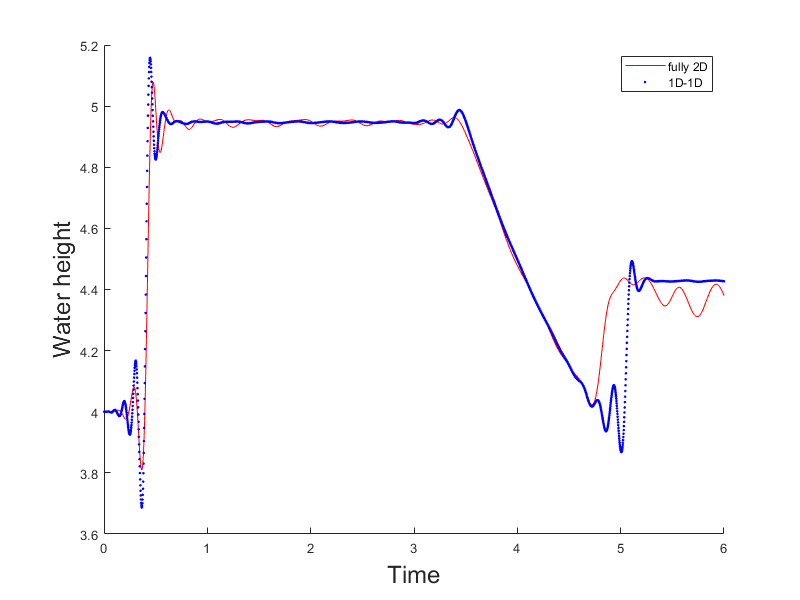}}
\hspace{.1em}
\subfloat[P3]{\includegraphics[width=.32\textwidth]{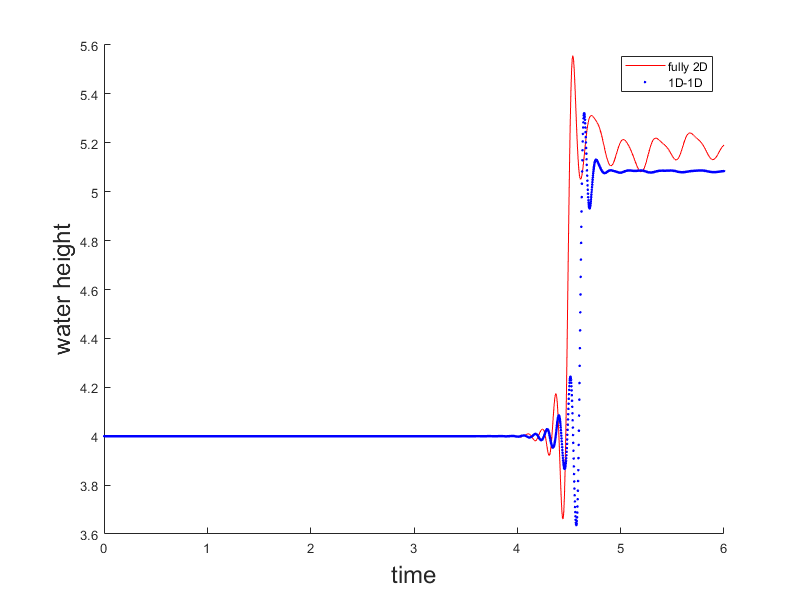}}
\caption{Water height in T-junction with initial conditions (\ref{eq:SWE_T_initial_condition_2}) at points P1 (a), P2 (b), and P3 (c) for the fully 2D and 1D-1D junction models.}
\label{fig:SWE_T_points_2}
\end{figure}
\begin{figure}
\centering
\subfloat[P1]{\includegraphics[width=.32\textwidth]{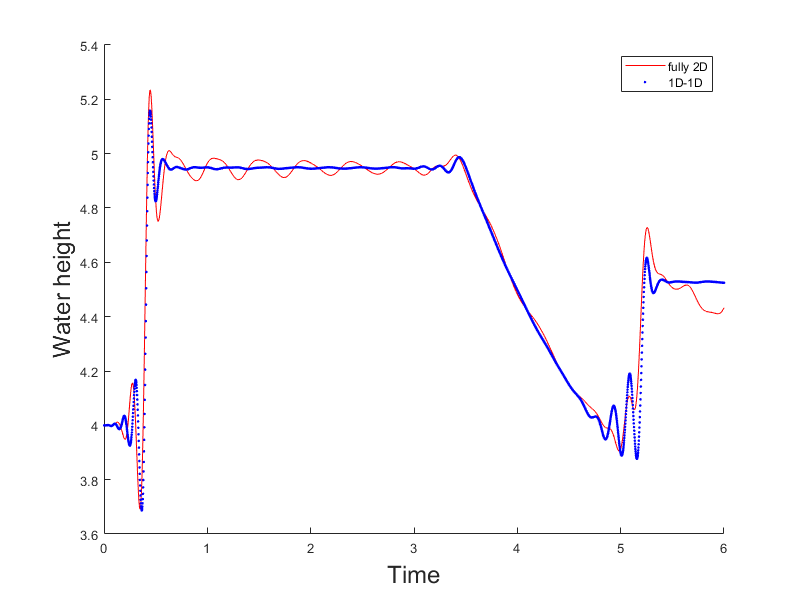}}
\hspace{.1em}
\subfloat[P2]{\includegraphics[width=.32\textwidth]{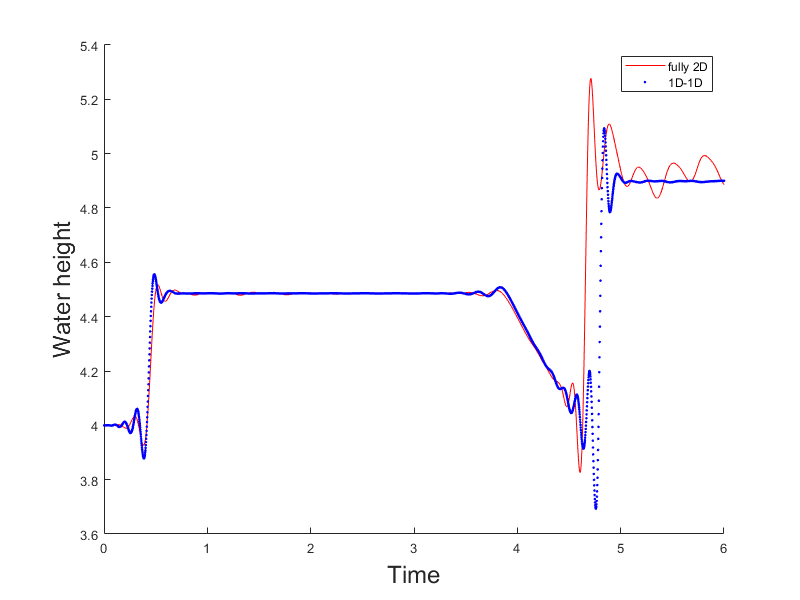}}
\hspace{.1em}
\subfloat[P3]{\includegraphics[width=.32\textwidth]{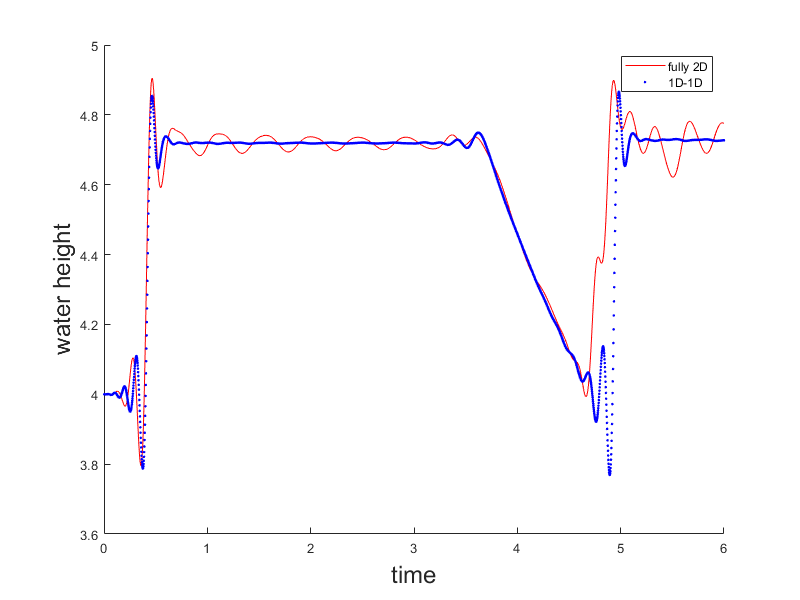}}
\caption{Water height in T-junction with initial conditions (\ref{eq:SWE_T_initial_condition_3}) at points P1 (a), P2 (b), and P3 (c) for the fully 2D and 1D-1D junction models.}
\label{fig:SWE_T_points_3}
\end{figure}

We also test with a smooth initial conditions (\ref{eq:SWE_T_initial_condition_4}) where we perturb a constant water height with a small sine wave. We plot the water height at $P1$, $P2$, and $P3$ in Figure \ref{fig:SWE_T_points_4}.
\begin{align}
h_0 = \begin{cases}
      4+0.1\sin(\pi x)  & \text{in channel 1}\\
      4+0.1\sin(\pi y)  & \text{in channel 2}\\
      4+0.1\sin(-\pi y)  & \text{in channel 3}\\
    \end{cases}, 
\hspace{1cm} u_0=v_0=0.
\label{eq:SWE_T_initial_condition_4}
\end{align}
\begin{figure}
\centering
\subfloat[P1]{\includegraphics[width=.32\textwidth]{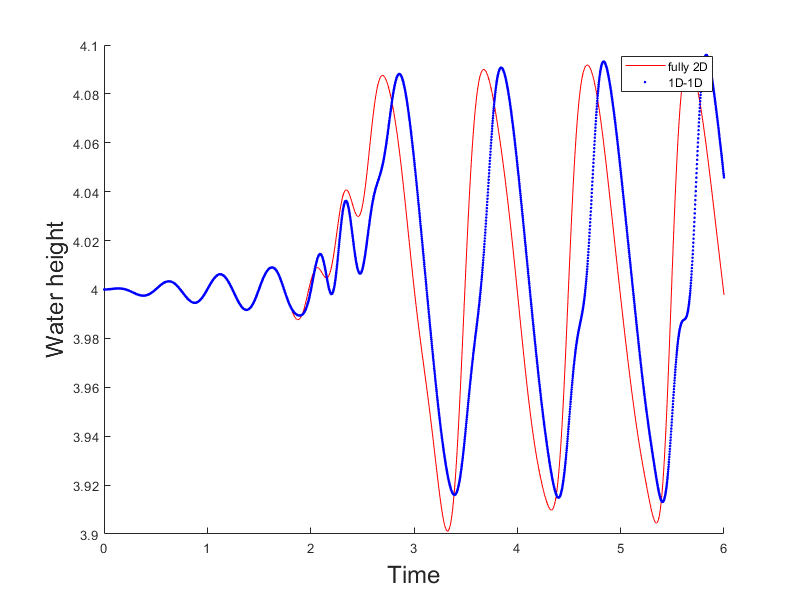}}
\hspace{.1em}
\subfloat[P2]{\includegraphics[width=.32\textwidth]{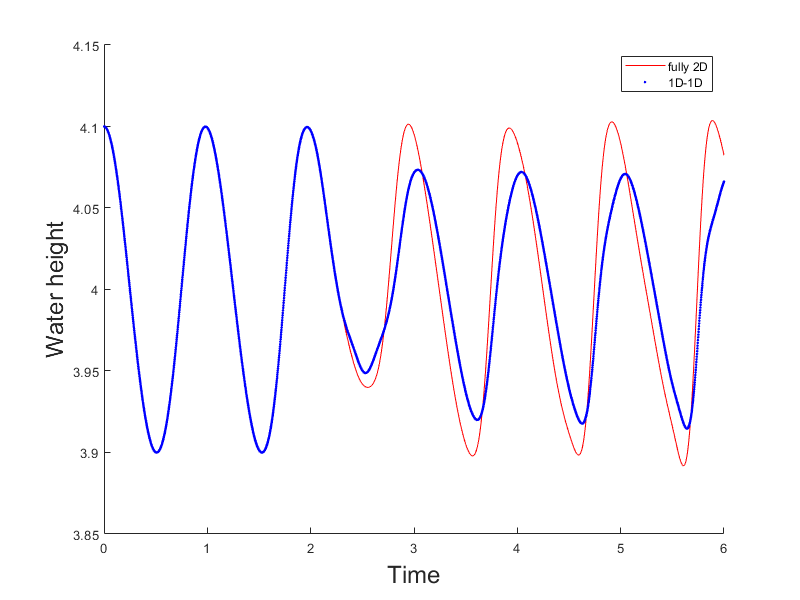}}
\hspace{.1em}
\subfloat[P3]{\includegraphics[width=.32\textwidth]{figs/swe_T_P2_4.png}}
\caption{Water height in T-junction with initial conditions (\ref{eq:SWE_T_initial_condition_4}) at points P1 (a), P2 (b), and P3 (c) for the fully 2D and 1D-1D junction models.}
\label{fig:SWE_T_points_4}
\end{figure}

We notice that in Figure \ref{fig:SWE_T_points_1}, \ref{fig:SWE_T_points_2}, and \ref{fig:SWE_T_points_3}, the water heights under the 1D-1D junction models approach constant values towards the end of the simulation, while the solutions from fully 2D model display more oscillatory behavior. From other numerical experiments under the fully 2D model, we observe that the shock wave from channel 1 hits the wall on the right side boundary and reflects back. The reflected wave, however, does not simply travel in the opposite direction, but spreads out in semicircular shape. The 1D model cannot capture this 2D behavior and thus fails to accurately predict the reflected wave. We also observe discrepancies in amplitude when using the 1D-1D junction treatment for the continuous initial conditions (\ref{eq:SWE_T_initial_condition_4}) in Figure \ref{fig:SWE_T_points_4}.

\subsubsection{Turning channel (1D-1D junction treatment)}
\begin{figure}[H]
\begin{center}
\includegraphics[width=.44\textwidth]{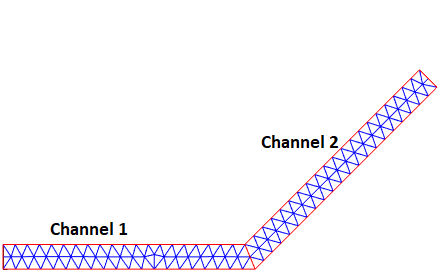}
\includegraphics[width=.28\textwidth]{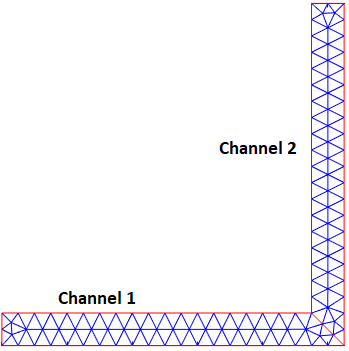}
\caption{2D mesh of $45^{\circ}$ and $90^{\circ}$ degree turning channel.}
\label{fig:2D_mesh_turns}
\end{center}
\end{figure}
We continue testing 1D-1D junction treatments by investigating the effect of channel angle for turning channels. In the turning channel experiment, we use a channel with a $45^{\circ}$ and a channel with a $90^{\circ}$ turn, as shown in Figure \ref{fig:2D_mesh_turns}. Each channel has a width of 1 and length of 10. We choose $P1$ and $P1$ to be the midpoint of each channel. PDEModel is used to construct the 2D mesh, with $h_{min} = 0.25$. For 1D-1D model, we represent each channel with 32 uniform elements of size 0.3125. We use polynomial degree $N=3$ in all experiments and run the experiment up to time $T=6$. We turn on Lax-Friedrichs dissipation to compare solutions at midpoints of the each 1D channel with solutions from the fully 2D model. We test on the same initial conditions for the shallow water equations:

\begin{align}
h_0 = \begin{cases}
      6  & \text{in $x\leq4$ channel 1}\\
      4  & \text{otherwise}\\
    \end{cases}, 
\hspace{1cm} u_0=v_0=0.
\label{eq:SWE_turn_initial_condition_1}
\end{align}
These initial conditions contain discontinuities. We plot the water heights at $P1$ and $P2$ in Figure \ref{fig:SWE_turn45_points}, \ref{fig:SWE_turn90_points}.
\begin{figure}[H]
\begin{center}
\subfloat[P1]{\includegraphics[width=.45\textwidth]{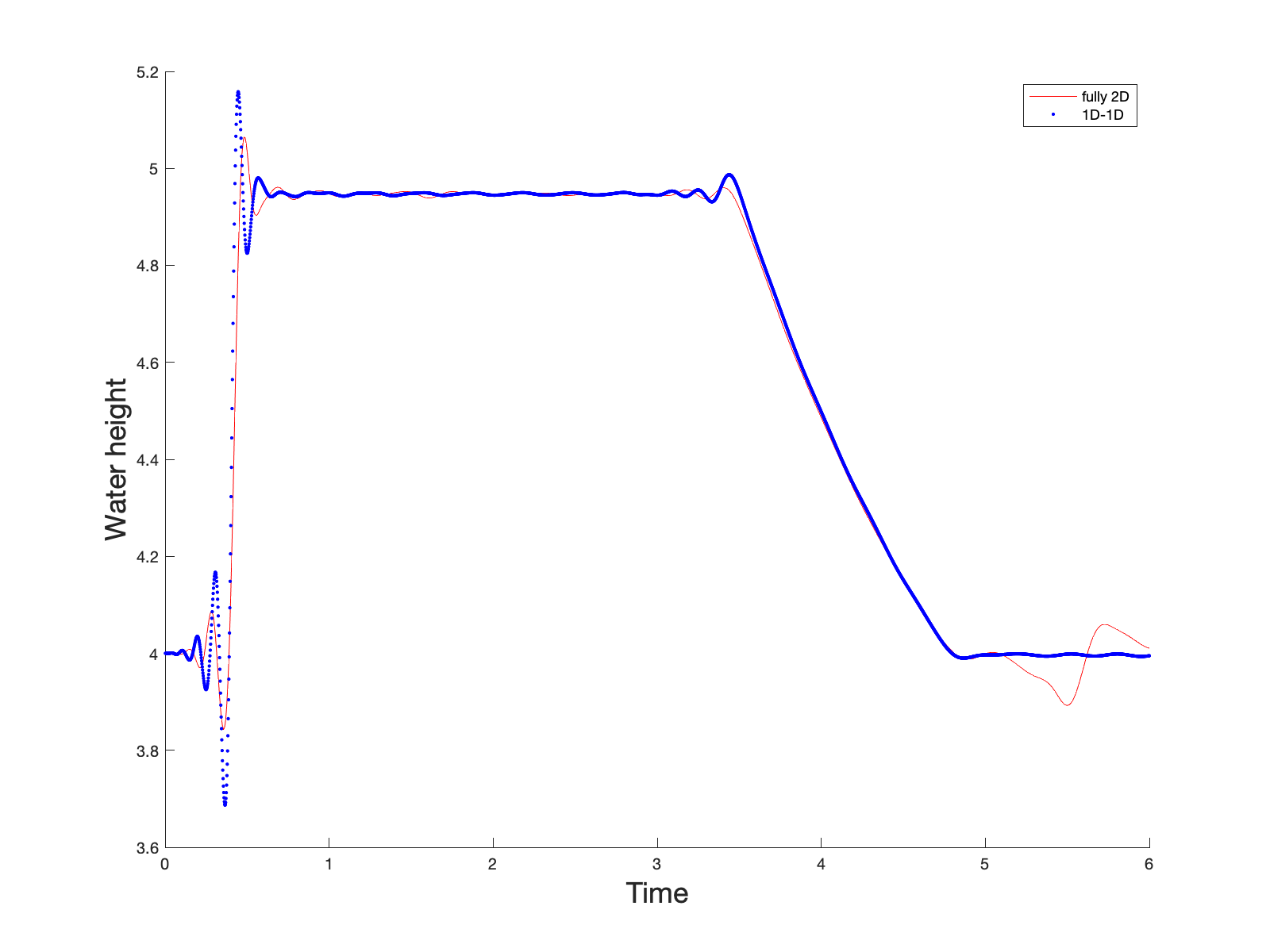}}
\subfloat[P2]{\includegraphics[width=.45\textwidth]{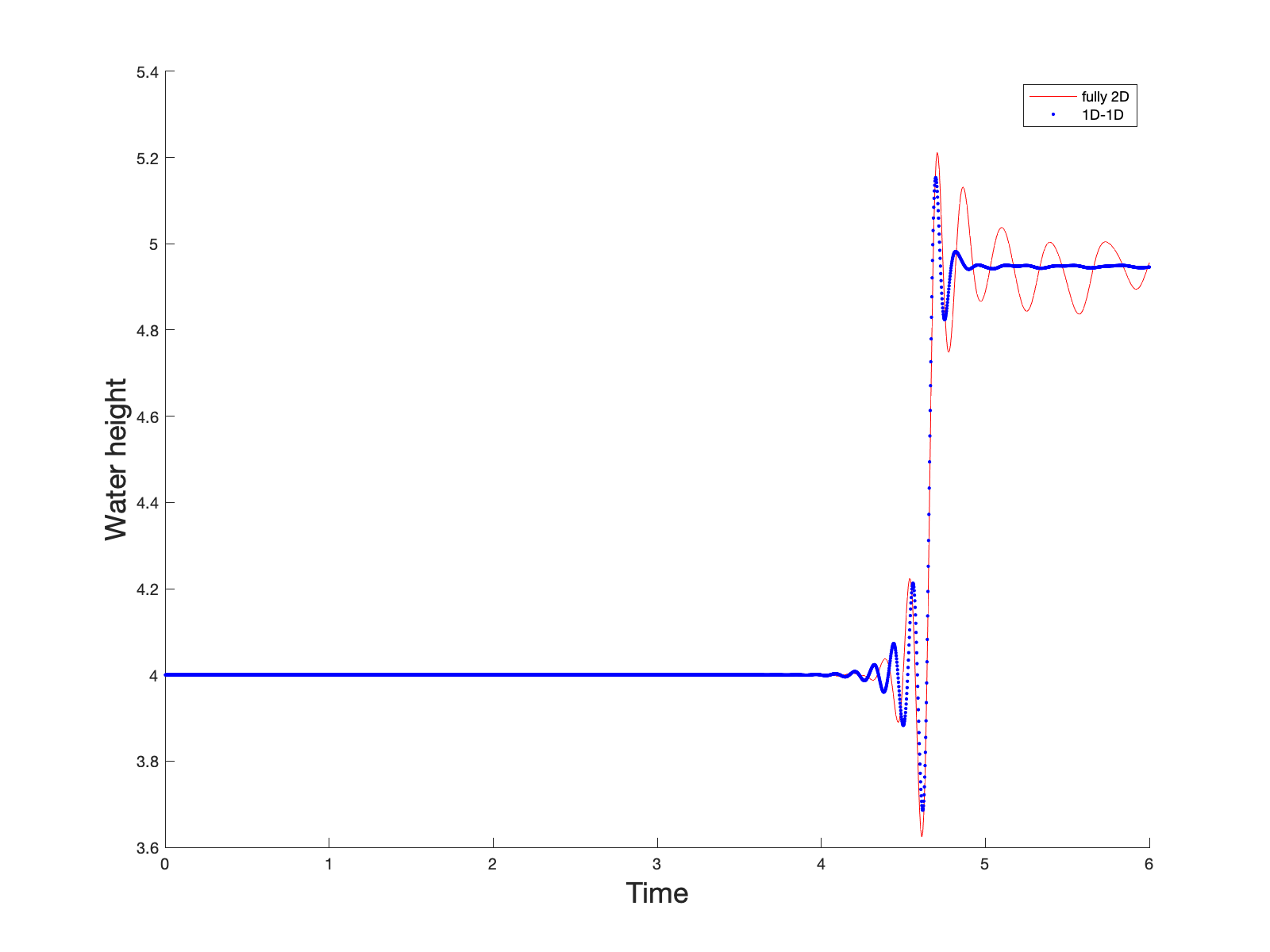}}
\caption{Water height in $45^{\circ}$ turn at points P1 (a) and P2 (b) for the fully 2D and 1D-1D junction models.}
\label{fig:SWE_turn45_points}
\end{center}
\end{figure}

\begin{figure}[H]
\begin{center}
\subfloat[P1]{\includegraphics[width=.45\textwidth]{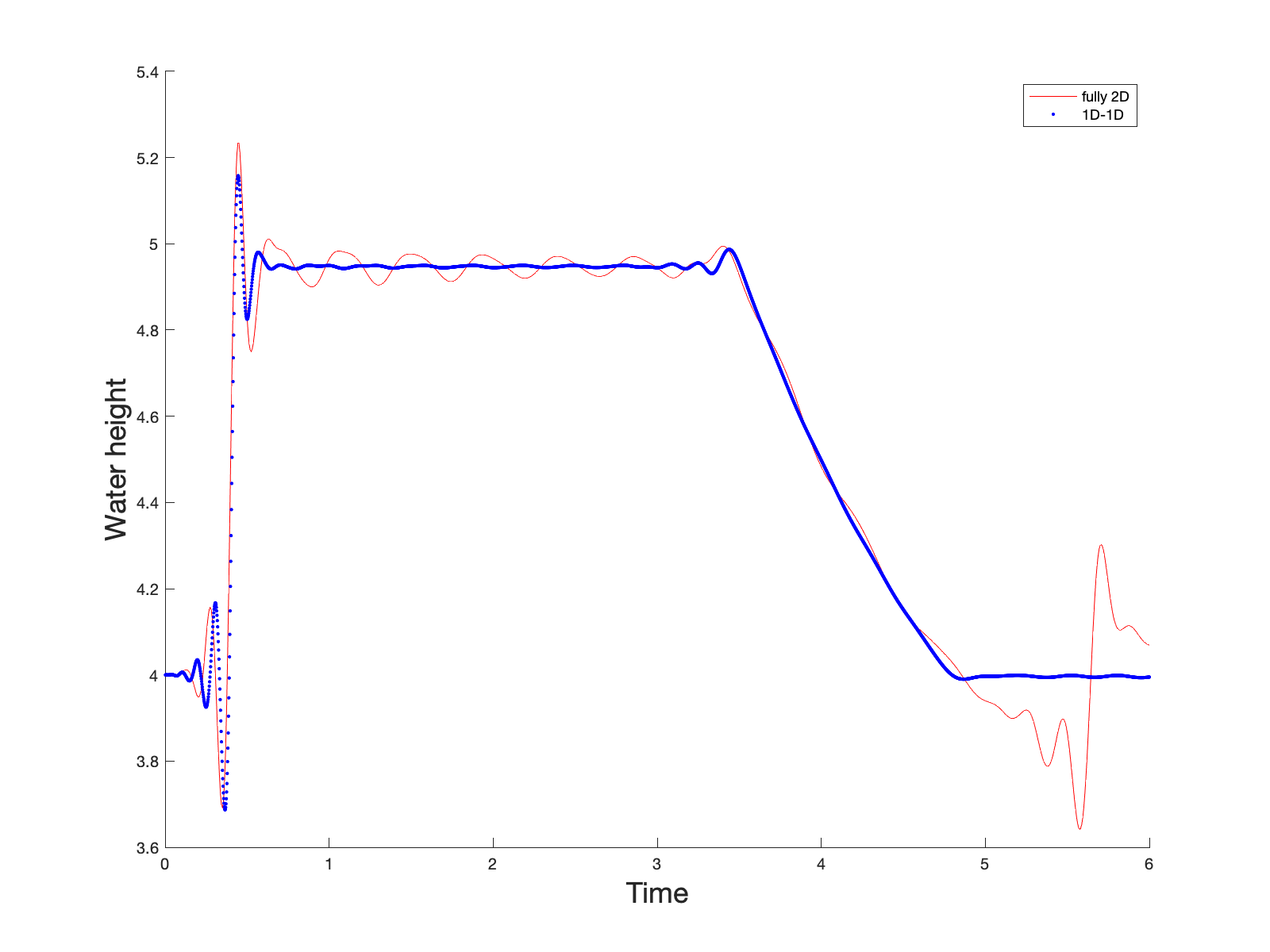}}
\subfloat[P2]{\includegraphics[width=.45\textwidth]{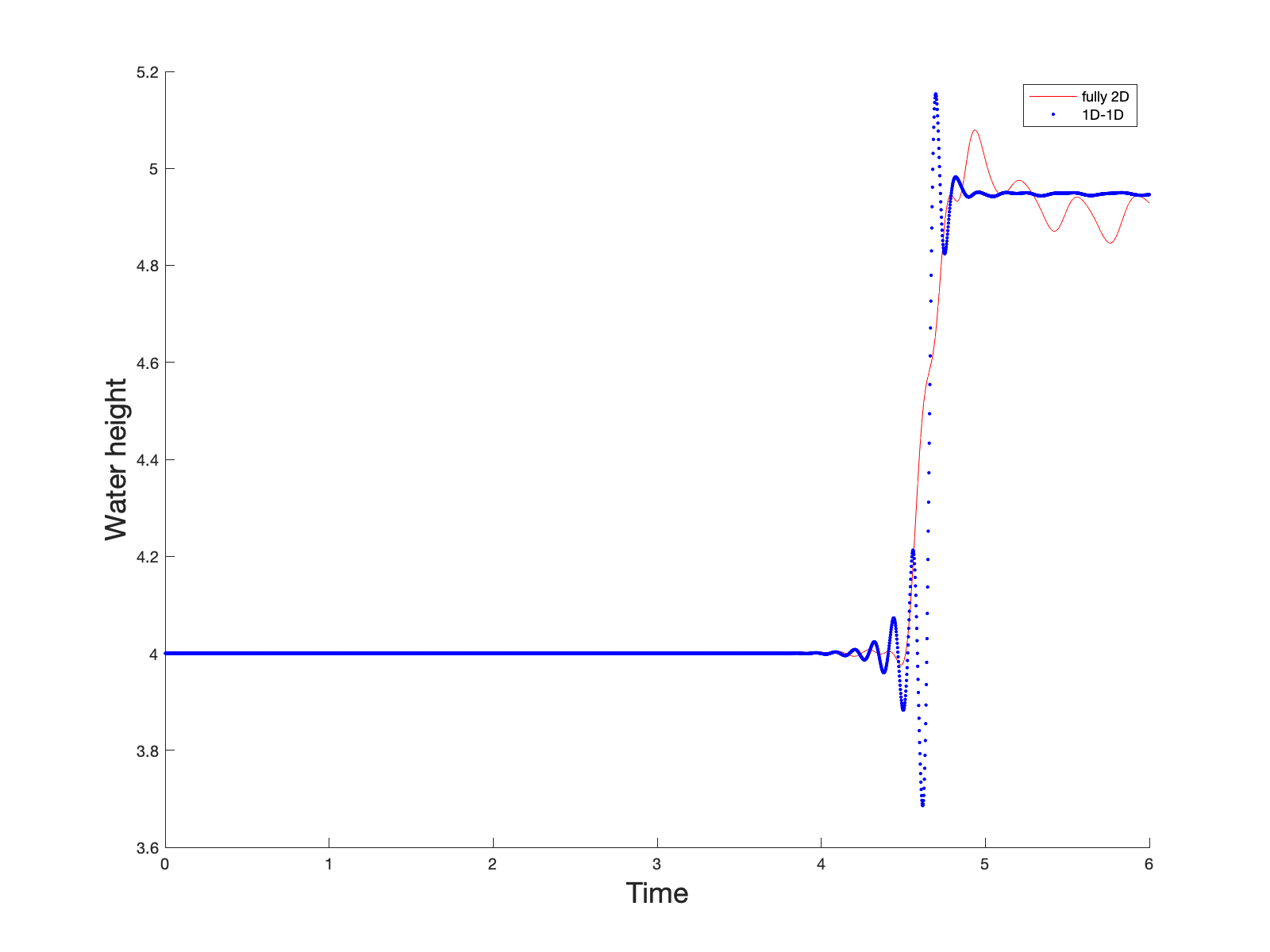}}
\caption{Water height in $90^{\circ}$ turn at points P1 (a) and P2 (b) for the fully 2D and 1D-1D junction models.}
\label{fig:SWE_turn90_points}
\end{center}
\end{figure}

We observe that the turning channel reflects part of the water back into channel 1. The magnitude of the reflected wave is greater for a larger turning angle. The 1D model approximates the fully 2D model reasonably well, but does not capture the reflected wave. We note that it may be possible to model this behavior using a partial wall boundary condition. However, this is beyond the scope of the paper and will be explored in the future work.

\subsubsection{Dam break and turning channel (1D-2D junction treatment)}
\label{subsec:dam_swe}
In this numerical experiment, we model a dam break and turning channel from \cite{bellamoli2018numerical} as shown in Figure \ref{fig:Dam_break_mesh}. The reservoir is represented by a square with dimension $2.5 \times 2.5$ and the channel (with a width of $0.5$) is connected to the right side of the reservoir. The length of each channel is 4. We also use the mesh generator within PDEModel to construct the 2D mesh, with $h_{min} = 0.25$. For the 1D-2D junction treatment, we model the reservoir using a 2D mesh. The point at which the channel turns is also modeled as a 2D domain. In the 1D mesh, we have 16 uniform elements for each channel of size 0.25. Wall boundary conditions are imposed on all boundaries for each simulation. We use polynomial degree $N=3$ in all experiments. The initial conditions for the shallow water equations are taken to be the following:
\begin{align}
h_0 = \begin{cases}
      10 & \text{in reservoir}\\
      6  & \text{in channel}\\
    \end{cases}, 
\hspace{1cm} u_0=v_0=0.
\label{eq:SWE_dam_initial_condition_1}
\end{align}
\begin{figure}[H]
\centering
\subfloat[]{\raisebox{8ex}{\includegraphics[height=.1\textheight]{figs/dam_channel_mesh.png}}}
\subfloat[]{\raisebox{8ex}{\includegraphics[height=0.0875\textheight]{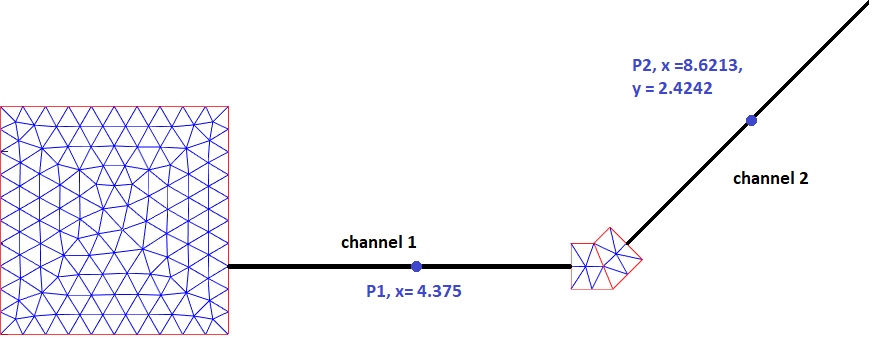}}}
\subfloat[]{\includegraphics[width=.32\textwidth]{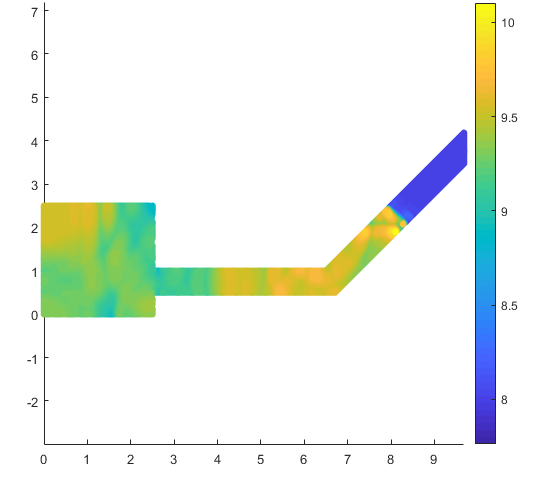}}
\caption{Dam break and turning channel meshes with initial conditions (\ref{eq:SWE_dam_initial_condition_1}). Fully 2D mesh (a), 1D-2D (b) and solution snapshot (c).}
\label{fig:Dam_break_mesh}
\end{figure}
We verify that the entropy RHS (\ref{eq:entropy_rhs}) is on the order of $10^{-15}$ to $10^{-13}$ on absence of Lax-Friedrichs dissipation. Then, we turn on the Lax-Friedrichs dissipation to compare the solutions at midpoints of the two 1D channels as marked in Figure \ref{fig:Dam_break_mesh}. We run the model up to time $T = 2.5$ and plot the water height over time in Figure \ref{fig:Dam_break_swe_P}. We observe that solutions from the fully 2D and 1D-2D junction models are similar.
\begin{figure}
\begin{center}

\subfloat[P1]{\includegraphics[width=.45\textwidth]{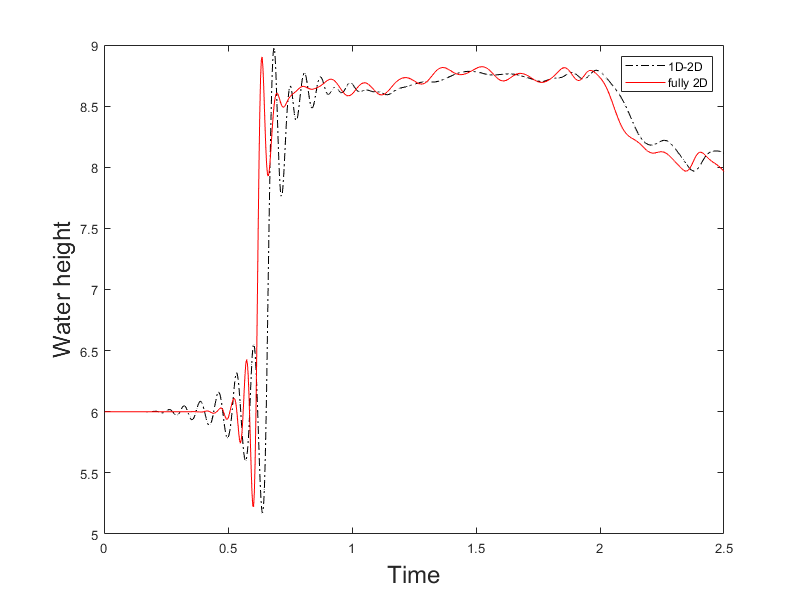}}
\subfloat[P2]{\includegraphics[width=.45\textwidth]{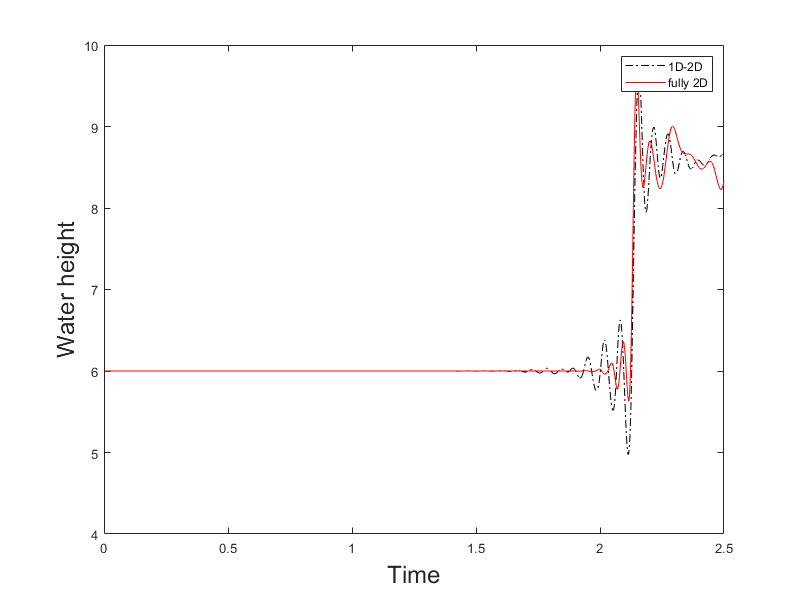}}

\caption{Water height in dam break and turning channel with continuous initial conditions (\ref{eq:SWE_dam_initial_condition_1}) at points P1 (a) and P2 (b) for the fully 2D and 1D-2D junction models.}
\label{fig:Dam_break_swe_P}
\end{center}
\end{figure}

We also test the shallow water equations with the following continuous initial conditions:
\begin{align}
h_0 = \begin{cases}
      2+0.1\sin(2\pi(x-4.375)/3.75) & \text{in channel 1}\\
      2  & \text{otherwise}\\
    \end{cases}, 
\hspace{1cm} u_0=v_0=0.
\label{eq:SWE_dam_initial_condition_2}
\end{align}
We then plot the results using this continuous initial conditions (\ref{eq:SWE_dam_initial_condition_2}) at $P1$ and $P2$ in Figure \ref{fig:Dam_break_swe_P_c}. Again, the fully 2D and 1D-2D junction models produce comparable results.
\begin{figure}[H]
\begin{center}
\subfloat[P1]{\includegraphics[width=.45\textwidth]{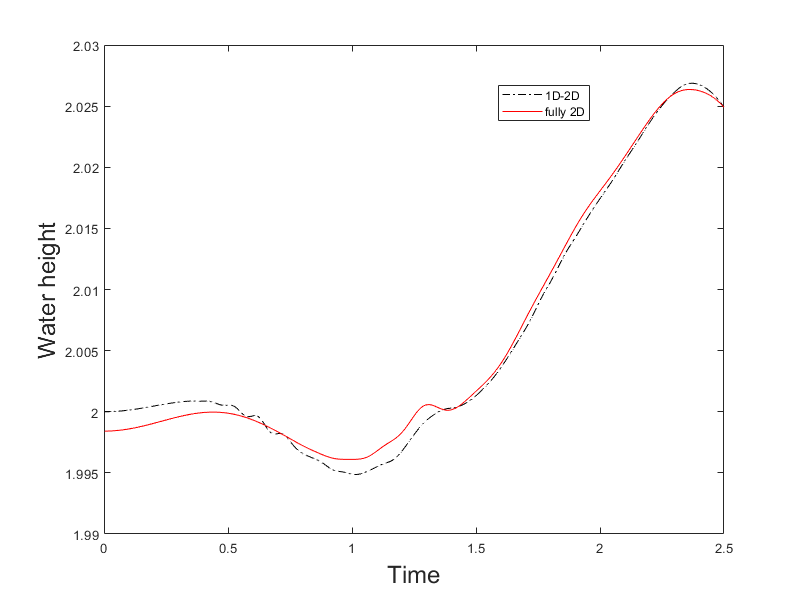}}
\subfloat[P2]{\includegraphics[width=.45\textwidth]{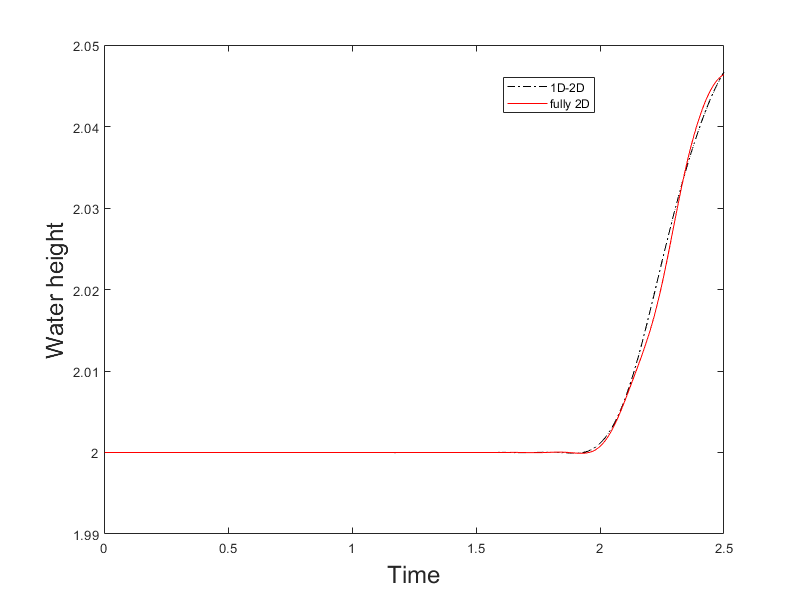}}
\caption{Water height in dam break and turning channel with discontinuous initial conditions (\ref{eq:SWE_dam_initial_condition_2}) at points P1 (a) and P2 (b) for the fully 2D and 1D-2D junction models.}
\label{fig:Dam_break_swe_P_c}
\end{center}
\end{figure}

\end{subequations}

%% file: section6.tex
\section{Conclusions}
\label{sec:Conclusion}
\hspace{1em}In this work, we present high-order entropy stable DG schemes for coupling 1D channels together. We construct both 1D-1D and 1D-2D junction treatments. We prove conservation of entropy for both junction treatments and describe how to apply entropy dissipation using Lax–Friedrichs interface penalization terms. We verify these results with numerical experiments for the shallow water and compressible Euler equations. We also compare numerical results on several different geometries with different initial conditions and conclude that 1D-2D model produces results similar to the fully 2D model, matching the observations in \cite{bellamoli2018numerical, neupane2015discontinuous}, whereas, the performance of the 1D-1D model is more sensitive to the domain geometry and initial conditions. This is due to the fact that the 1D-1D model cannot capture higher dimensional flows behavior (e.g., flow in the direction perpendicular to each 1D channel). Apart from computational cost, the only advantage of the 1D-1D junction is simplicity of implementation. The 1D-2D junction model is more robust to higher-dimensional flow effects at junctions; however, 1D-1D junction models can be used for flows which are known to be less sensitive to junction geometry. We note that this paper does not address the issue related to the accuracy of 1D-1D and 1D-2D junction models, but on ensuring that such junctions can be modeled numerically in an entropy stable fashion. Since comparisons of 1D-1D and 1D-2D junction models lie outside the scope of this paper, we refer interested readers to \cite{neupane2015discontinuous, bellamoli2018numerical} for comparisons between types of junction models.

%% file: appendixA.tex
\appendix
\section*{Appendix A}
\label{sec:Euler_experiments}
\renewcommand{\thesubsection}{\Alph{subsection}}
\renewcommand{\theequation}{A.\arabic{equation}}
In this appendix, we examine junction treatments for the compressible Euler equations with entropy conservative fluxes. The compressible Euler equations for gas dynamics in two dimensions are given by
\begin{align}
\frac{\partial }{\partial t}\begin{bmatrix}
\rho\\
\rho u\\
\rho v\\
E
\end{bmatrix}+
\frac{\partial }{\partial x}\begin{bmatrix}
\rho u\\
\rho u^2+p\\
\rho uv\\
u(E+p)
\end{bmatrix}+
\frac{\partial }{\partial y}\begin{bmatrix}
\rho v\\
\rho uv\\
\rho v^2+p\\
v(E+p)
\end{bmatrix} = 0.
\label{eq:Euler_2D}
\end{align}
Here, $\rho$ and $p$ denote density and the pressure, respectively. The velocity in the $x$ direction is denoted by $u$ and the velocity in the $y$ direction is denoted by $v$. The total energy is denoted by $E$ and satisfies the constitutive relation involving the pressure $p$
\begin{align}
E = \frac{1}{2}\rho\nor{U}^2 + \frac{p}{\gamma-1},
\end{align}
where $\nor{U}^2 = u^2+v^2$, and $\gamma = 1.4$ is the ratio of specific heat a diatomic gas. In this example, we have conservative variables $\bm{u} = [\rho, \rho u, \rho v,E]^T$ and flux functions $\bm{f}_1 = [\rho u, \rho u^2+p, \rho uv, u(E+p)]^T$ and $\bm{f}_2 = [\rho v, \rho uv, \rho v^2+p, v(E+p)]^T$.

The one-dimensional compressible Euler equations can also be derived under assumptions similar to those used to derive the one-dimensional shallow water equations from the two-dimensional system. In one dimension, the compressible Euler equations are
\begin{align}
\frac{\partial }{\partial t}\begin{bmatrix}
\rho\\
\rho u\\
E
\end{bmatrix}+
\frac{\partial }{\partial x}\begin{bmatrix}
\rho u\\
\rho u^2+p\\
u(E+p)
\end{bmatrix}= 0.
\label{eq:Euler_1D}
\end{align}
where we define $\nor{U}^2 = u^2$  in one dimension.

The transform matrix $\bm{R}$ between 1D and 2D for the compressible Euler  equations is
\begin{align}
\bm{R} = \begin{bmatrix}
1 &0 &0\\
0 &n_1 &0\\
0 &n_2 &0\\
0 &0 &1\\
\end{bmatrix}, \qquad \bm{R}^T\bm{R} = \bm{I}.
\label{eq:1d2d_transformation_matrix_euler}
\end{align}

In this work, the mathematical entropy for the compressible Euler equations is taken to be the unique mathematical entropy for the compressible Navier-Stokes equations \cite{hughes1986new}
\[
S(\bm{u}) = -\rho s, 
\]
where $s = \log\LRp{\frac{p}{\rho^\gamma}}$ is the physical specific entropy. The entropy variables $\bm{v}$ in $d$ dimensions are 
\begin{align}
v_1 = \frac{\rho e (\gamma + 1 - s) - E}{\rho e}, \qquad v_{1+ i}= \frac{\rho {{u}_i}}{\rho e}, \qquad v_{d+2} = -\frac{\rho}{\rho e}, \label{eq:evars}
\end{align}
for $i = 1,\ldots, d$.  The inverse map from entropy to conservative variables is 
\begin{align*}
\rho = -(\rho e) v_{d+2}, \qquad 
\rho {u_i} = (\rho e) v_{1+i}, \qquad 
 E = (\rho e)\LRp{1 - \frac{\sum_{j=1}^d{v_{1+j}^2}}{2 v_{d+2}}},
\end{align*}
where $i = 1,\ldots,d$, and $\rho e$ and $s$ in terms of the entropy variables are 
\begin{equation*}
\rho e = \LRp{\frac{(\gamma-1)}{\LRp{-v_{d+2}}^{\gamma}}}^{1/(\gamma-1)}e^{\frac{-s}{\gamma-1}}, \qquad 
s = \gamma - v_1 + \frac{\sum_{j=1}^d{v_{1+j}^2}}{2v_{d+2}}.
\end{equation*}

To introduce the entropy conservative fluxes for the compressible Euler equations, we start with some notations. Let $f$ denote some piecewise continuous function, and $f^+$ denote the exterior value of $f$ across an element face. The average and logarithmic averages are
\begin{align}
\avg{f} = \frac{f+f^+}{2}, \hspace{1cm} \avg{f}^{\rm{log}} = \frac{f^+-f}{\rm{log}(f^+)-\rm{log} (f)}.
\end{align}
The average and logarithmic average are assumed to act component-wise on vector valued functions.

The entropy conservative numerical fluxes for the 2D compressible Euler equations are given by Chandrashekar \cite{chandrashekar2012kinetic}:
\begin{align}
\bm{f}^x_{S}\LRp{\bm{u}_L,\bm{u}_R} &=
\begin{bmatrix}
\avg{p}^{\rm{log}}\avg{u}\\
\avg{p}^{\rm{log}}\avg{u}^2 + p_{\rm{avg}}\\
\avg{p}^{\rm{log}}\avg{u}\avg{v}\\
(E_{\rm{avg}} + p_{\rm{avg}}) \avg{u}
\end{bmatrix}, 
\qquad
\bm{f}^y_{S}\LRp{\bm{u}_L,\bm{u}_R} &=
\begin{bmatrix}
\avg{p}^{\rm{log}}\avg{v}\\
\avg{p}^{\rm{log}}\avg{u}\avg{v}\\
\avg{p}^{\rm{log}}\avg{v}^2 + p_{\rm{avg}}\\
(E_{\rm{avg}} + p_{\rm{avg}}) \avg{v}
\end{bmatrix}.
\label{eq:Euler_flux2d}
\end{align}
where we need to introduce the auxiliary quantities $\beta = \frac{\rho}{2p}$ and 
\begin{align}
p_{\rm{avg}} = \frac{\avg{\rho}}{2 \avg{\beta}} , \hspace{1cm} E_{\rm{avg}} =
\frac{\avg{\rho}^{\rm{log}}}{2 \avg{\beta}^{\rm{log}} (\gamma-1)} + \frac{u_{\rm{avg}}^2}{2}, \hspace{1cm} u_{\rm{avg}}^2 = u_Lu_R+v_Lv_R.
\end{align}

The entropy conservative fluxes for the compressible Euler equations in 1D are
\begin{align}
\bm{f}_{S1D}\LRp{\bm{u}_L,\bm{u}_R} &=
\begin{bmatrix}
\avg{p}^{\rm{log}}\avg{u}\\
\avg{p}^{\rm{log}}\avg{u}^2 + p_{\rm{avg}}\\
(E_{\rm{avg}} + p_{\rm{avg}}) \avg{u}
\end{bmatrix},
\label{eq:Euler_flux1d}
\end{align}
where we need calculate the term $E_{\rm{avg}}$ with $u_{\rm{avg}}^2 = u_Lu_R$ in 1D
\subsection{Numerical experiments for the compressible Euler equations}
\subsubsection{Parallel split and converge (1D-2D and 1D-1D junction treatments)}
For the compressible Euler equations, we reuse the same mesh and setup in Section \ref{subsec:ps_swe} (as shown in Figure \ref{fig:2D_mesh_split}) with the following initial conditions:
\begin{align}
\rho_0 = \sin(\pi x/2)+2, \hspace{0.8cm} u_0=2,\hspace{0.8cm} v_0=0,  \hspace{0.8cm} \gamma  = 1.4, \hspace{0.8cm}
p_0 = 2.
\label{eq:euler_ps_init_cd_1}
\end{align}
We also test with different polynomial degree on each domain and list the maximum absolute value of the entropy RHS (\ref{eq:entropy_rhs}) up to time $T=1$. Results for the 1D-2D model are shown in Table \ref{Tab:Euler_entopy_rhs} and results for the 1D-1D model are shown in \ref{Tab:1D_1D_entopy_rhs}.
\begin{table}[H]
\begin{center}
\begin{tabular}{|c|c|c|c|}
\hline
                 &  $N_{2D} = 3$  &  $N_{2D} = 4$  &  $N_{2D} = 5$ \\
\hline
  $N_{1D} = 3$  &  2.0872e-14 & 4.9950e-14 & 1.3084e-13 \\
\hline
  $N_{1D} = 4$  &  1.1824e-13 & 2.3324e-13 & 3.6135e-13 \\
\hline
  $N_{1D} = 5$  &  3.6599e-13 & 4.5905e-13 & 7.9016e-13 \\
\hline
\end{tabular}
\end{center}
\caption{Maximum of absolute value of entropy RHS (\ref{eq:entropy_rhs}) for compressible Euler  1D-2D coupling.}
\label{Tab:Euler_entopy_rhs}
\end{table}

\begin{table}[H]
\begin{center}
\begin{tabular}{|c|c|c|c|}
\hline
                 &  $N_{1D} = 3$  &  $N_{1D} = 4$  &  $N_{1D} = 5$ \\
\hline
  SWE    &  1.1191e-13 & 7.5495e-14 & 8.3311e-13 \\
\hline
  Euler  &  2.3082e-13 & 1.7586e-13 & 2.9110e-13 \\
\hline
\end{tabular}
\end{center}
\caption{Maximum of absolute value of entropy RHS (\ref{eq:entropy_rhs}) for compressible Euler 1D-1D coupling.}
\label{Tab:1D_1D_entopy_rhs}
\end{table}
\begin{figure}[H]
\centering
\subfloat[P1]{\includegraphics[width=.49\textwidth]{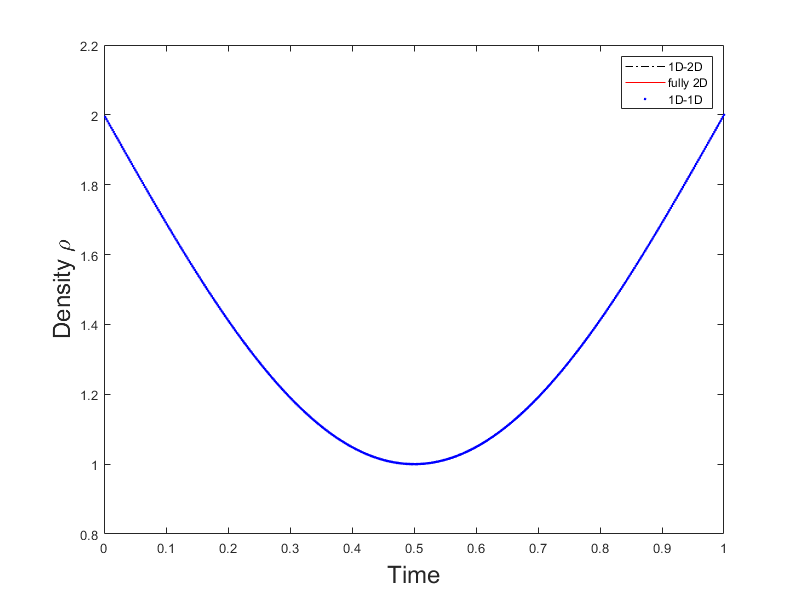}}
\subfloat[Error]{\includegraphics[width=.49\textwidth]{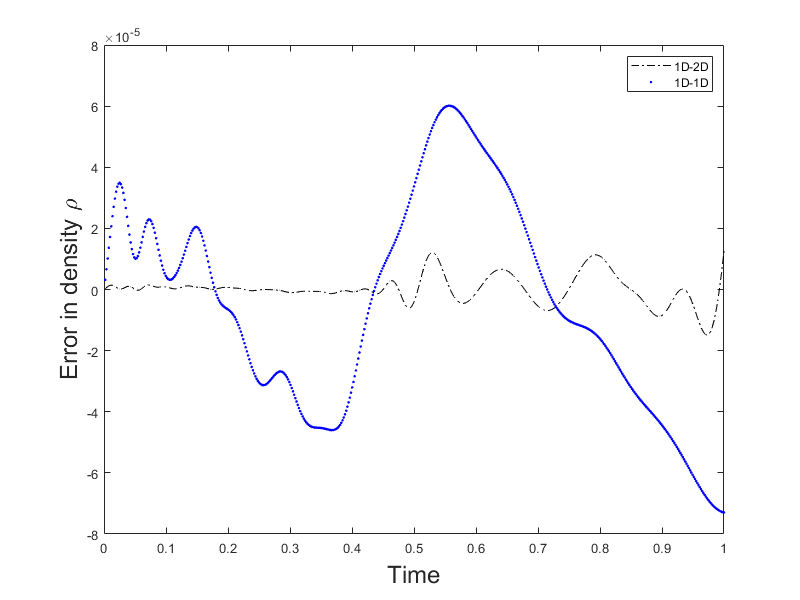}}
\caption{Density $\rho$ in the parallel split problem with initial conditions (\ref{eq:euler_ps_init_cd_1}) at point P1 for the fully 2D, 1D-2D, and 1D-1D junction models. Errors are computed using the fully 2D model as the ``exact'' solution.}
\label{fig:Euler_parallel_split_points_time}
\end{figure}
To test accuracy of these three models, fully 2D, 1D-2D and 1D-1D junction treatments, for the compressible Euler equations, we plot the results at  the midpoints of each channel, as marked in Figure \ref{fig:2D_mesh_split}. Coincidentally, because we have the periodic initial conditions and these points are separated by exactly one wavelength, all three points share the same solutions. We notice that with continuous solutions, all three models produce very solutions with small errors as shown in Figure \ref{fig:Euler_parallel_split_points_time}.

From these experiments, we can conclude that our numerical method is entropy conservative for both the shallow water equations and the compressible Euler equations using either 1D-2D or 1D-1D junction treatments. Different models produce different oscillations near the jump, but the their magnitudes are on the same scale. For the solutions that remain continuous, both 1D-2D and 1D-1D junction models generate solutions extremely close to the fully 2D model with absence of vertical flows. However, we expect the 1D-1D junction model to fail where fully 2D motions exist near the junction as in the shallow water experiment in Section \ref{sec:numerical_results}.

\subsubsection{Diamond split and converge}
For our second experiment with the compressible Euler equations, we reuse the diamond split setup in Section \ref{subsec:diamond_swe} (as shown in Figure \ref{fig:diamond_mesh}) with the following initial conditions:
\begin{align}
\rho_0 = 2, \hspace{0.8cm} u_0=v_0=0,  \hspace{0.8cm} \gamma  = 1.4, \hspace{0.8cm}
p_0 = \begin{cases}
      3 & \text{in channel 1}\\
      4 & \text{otherwise}\\
    \end{cases}.
\label{eq:Euler_diam_init_1}
\end{align}

\begin{figure}
\centering
\subfloat[]{\includegraphics[width=.49\textwidth]{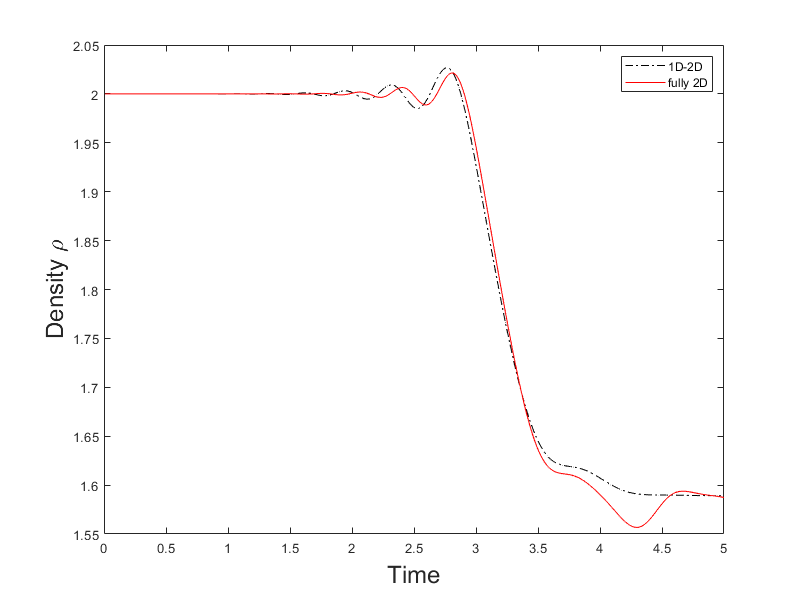}}
\subfloat[]{\includegraphics[width=.49\textwidth]{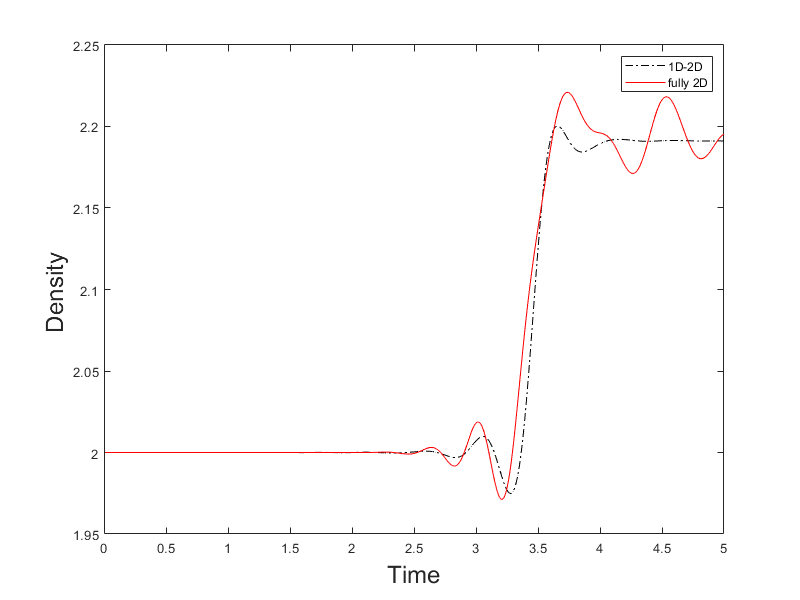}}
\caption{Density $\rho$ in diamond split and converge with initial conditions (\ref{eq:Euler_diam_init_1}) at points P1 (a) and P2 (b) for the fully 2D model and 1D-2D junction models.}
\label{fig:diamond_euler_p}
\end{figure}

We run the test without local Lax-Friedrichs penalization up to time $T=5$ to verify the conservation of entropy. Then, we enable local Lax-Friedrichs penalization for our accuracy test. We plot the values of $\rho$ at midpoints $P1$ and $P2$ from both fully 2D and 1D-2D junction treatments in Figure \ref{fig:diamond_euler_p}. 

We also test the compressible Euler equations with continuous initial conditions:
\begin{align}
\rho_0 = 2, \hspace{0.5cm} u_0=v_0=0,  \hspace{0.5cm} \gamma  = 1.4, \hspace{0.5cm}
p_0 = \begin{cases}
      2+\sin(\pi(x+5.5\sqrt{2}+5)/5) & \text{in channel 1}\\
      2 & \text{otherwise}\\
    \end{cases}.
\label{eq:Euler_diam_init_2}
\end{align}

\begin{figure}
\centering
\subfloat[]{\includegraphics[width=.49\textwidth]{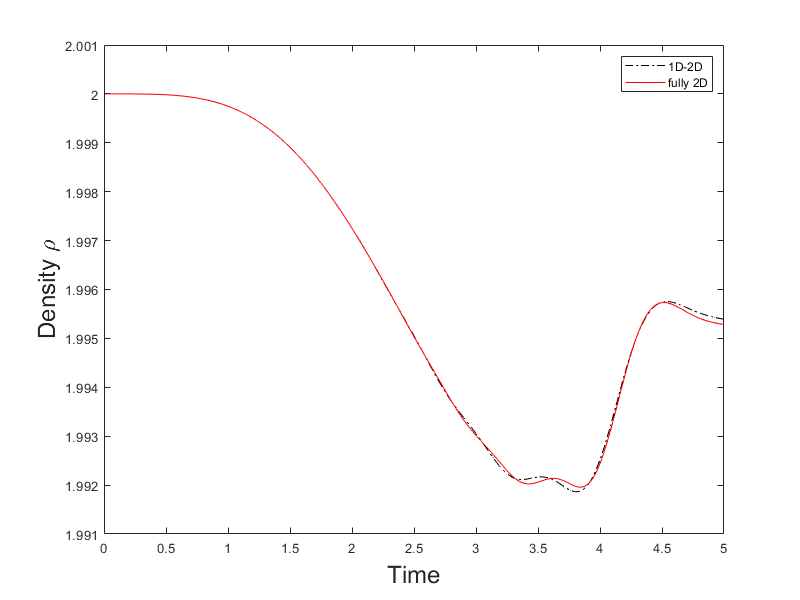}}
\subfloat[]{\includegraphics[width=.49\textwidth]{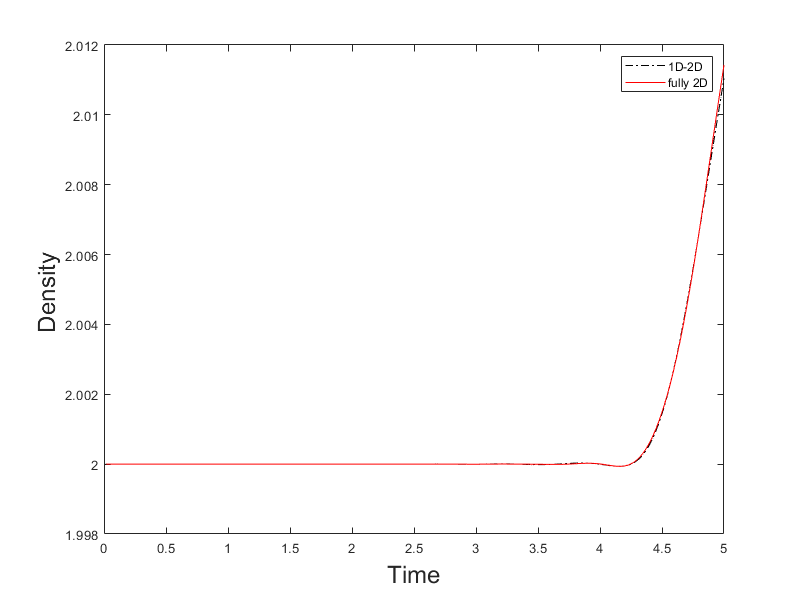}}
\caption{Density $\rho$ in the diamond split and converge with initial conditions (\ref{eq:Euler_diam_init_2}) at points P1 (a) and P2 (b) for the fully 2D and 1D-2D junction models.}
\label{fig:diamond_euler_p_c}
\end{figure}

We run the test up to time $T=5$ and plot the values of $\rho$ at $P1$ and $P2$ from fully 2D and 1D-2D junction models in Figure \ref{fig:diamond_euler_p_c}. 
In both shallow water and compressible Euler equations, we observe that the 1D-2D capture the general trend of the flow, but produce slightly different oscillation patterns compared to the fully 2D model.

\subsubsection{Dam break and turning channel}
\begin{figure}[H]
\begin{center}
\subfloat[]{\includegraphics[width=.45\textwidth]{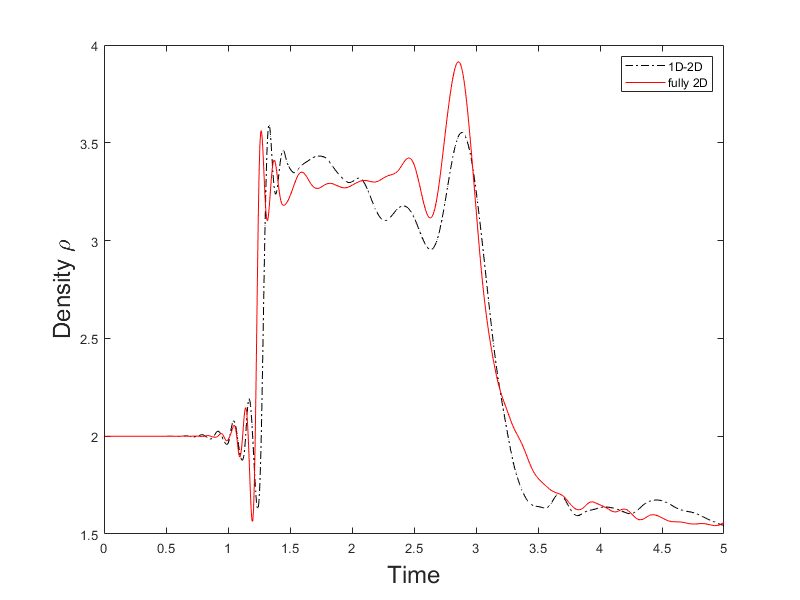}}
\subfloat[]{\includegraphics[width=.45\textwidth]{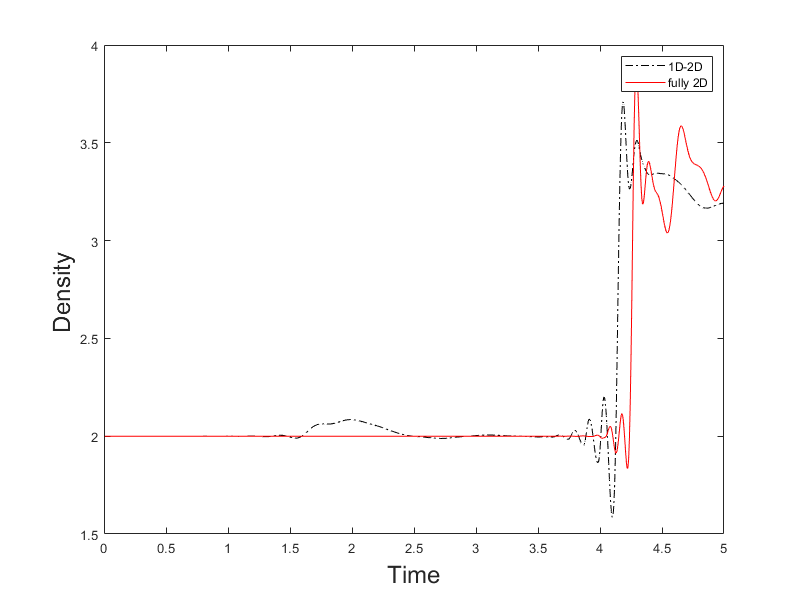}}
\caption{Density $\rho$ in the dam break and turning channel with initial conditions (\ref{eq:euler_dam_init}) at points P1 (a) and P2 (b) for the fully 2D and 1D-2D junction models.}
\label{fig:Dam_break_euler_P}
\end{center}
\end{figure}

Last, we test the compressible Euler equations on the dam break and turning channel setting in Section \ref{subsec:dam_swe}, as shown in Figure \ref{fig:Dam_break_mesh}. We first confirm conservation of entropy in the absence of Lax-Friedrichs penalization with the initial conditions:
\begin{align}
\rho_0 = 2, \hspace{0.8cm} u_0=v_0=0,  \hspace{0.8cm} \gamma  = 1.4, \hspace{0.8cm}
p_0 = \begin{cases}
      5 & \text{in reservoir }\\
      2 & \text{in channel}\\
    \end{cases}.
\label{eq:euler_dam_init}
\end{align}
We test the accuracy of the model with Lax-Friedrichs penalization and plot $\rho$ at midpoints $P1$ and $P2$ in the Figure \ref{fig:Dam_break_euler_P}. We run up to time $T=5$ and notice that similar patterns are generated from both the 1D-2D and fully 2D models. At $P1$, the oscillations from the fully 2D and 1D-2D junction models have noticeable discrepancies from time $T=1.5$ to $T=3.5$. We also find that there is a small bump at point $P2$ around time $T=2$, which the 1D-2D model does not capture.
